\documentclass[a4paper,11pt]{article}
\usepackage{pos}

\usepackage{tikz-cd}

\newtheorem{proposition}{Proposition}[section]
\newtheorem{definition}[proposition]{Definition}
\newtheorem{lemma}[proposition]{Lemma}
\newtheorem{theorem}[proposition]{Theorem}
\newtheorem{remark}[proposition]{Remark}
\newtheorem{example}[proposition]{Example}

\allowdisplaybreaks

\title{Gauge transformations on quantum principal bundles}

\author[a]{Antonio Del Donno}
\author[b,c]{Emanuele Latini}
\author*[a]{Thomas Weber}

\affiliation[a]{Mathematical Institute of Charles University,\\
 Sokolovsk\'a 83, Prague, Czech Republic}

\affiliation[b]{Dipartimento di Matematica, Università di Bologna,\\
Piazza di Porta S. Donato 5, Bologna, Italy}

\affiliation[c]{Istituto Nazionale di Fisica Nucleare, Sezione di Bologna,\\ 
via Berti Pichat 6/2, Bologna, Italy}

\emailAdd{antonio.deldonno@matfyz.cuni.cz}
\emailAdd{emanuele.latini@unibo.it}
\emailAdd{thomas.weber@matfyz.cuni.cz}

\abstract{We understand quantum principal bundle as faithfully flat Hopf--Galois extensions, with a structure Hopf algebra coacting on a total space algebra and with base algebra given by the coinvariant elements. To endow such bundles with a compatible differential structure, one requires the coaction to extend as a morphism of differential graded algebras. This leads to an exact noncommutative Atiyah sequence, a graded Hopf--Galois extension of differential forms and a canonical braiding on total space forms such that the latter are graded-braided commutative. We recall this approach to noncommutative differential geometry and further discuss the extension of quantum gauge transformations, in the sense of Brzezi\'nski, to differential forms. In this way we obtain an action of quantum gauge transformations on connections of the quantum principal bundle and their curvature. Explicit examples, such as the noncommutative 2-torus, the quantum Hopf fibration and smash product algebras are discussed.}

\FullConference{Proceedings of the Corfu Summer Institute 2024 "School and Workshops on Elementary Particle Physics and Gravity" (CORFU2024)\
12 - 26 May, and 25 August - 27 September, 2024\\
Corfu, Greece\\}

\tableofcontents

\begin{document}
\maketitle

\section{Introduction}

In the algebraic approach to noncommutative geometry, Hopf algebras play an essential role. As "quantum groups", they unify group algebras and universal enveloping algebras of Lie algebras and generalize them to the noncommutative realm. This sparked the interest of mathematicians and mathematical physicists alike, latest after the contribution of Drinfel'd \cite{Drinfeld}. In the seminal paper \cite{Woronowicz} Woronowicz introduced differential structures on Hopf algebras $H$ and characterized bicovariant differential calculi $\Omega^\bullet(H)$. As noticed by Schauenburg \cite{Schauenburg}, such covariant differential calculi admit a graded Hopf algebra structure. In a parallel development \cite{Schneider}, Schneider interpreted Hopf--Galois extensions $B\subseteq A$ as noncommutative generalizations of principal bundles, where a structure Hopf algebra $H$ coacts on a comodule algebra $A$ with coinvariant subalgebra $B=A^{\mathrm{co}H}$. The canonical Galois map encodes principality of such a bundle and faithful flatness replaces properness of a classical Lie group action. This gave a comprehensive geometrical picture to the algebraic approach. In their highly influential work \cite{BrzMaj}, Brzezi\'nski and Majid endowed Hopf--Galois extensions (or, more in general, comodule algebras), with first order differential calculi. This was crafted in a way that exactness of a noncommutative Atiyah sequence was ensured. An alternative, though closely related, approach to differential calculi on quantum principal bundles was developed by \DJ ur\dj evi\'c \cite{DurII}. He defined "complete" differential calculi $\Omega^\bullet(A)$ on Hopf--Galois extensions as differential calculi such that the coaction $\Delta_A\colon A\to A\otimes H$ extends to a (necessarily unique) differential graded algebra morphism $\Delta_A^\bullet\colon\Omega^\bullet(A)\to\Omega^\bullet(A)\otimes\Omega^\bullet(H)$. In this framework, vertical forms, horizontal forms and basic forms can be defined and completeness of the total space calculus automatically leads to exactness of the noncommutative Atiyah sequence
$$
0\to\mathrm{hor}^1\hookrightarrow\Omega^1(A)\overset{\pi_v}{\twoheadrightarrow}\mathrm{ver}^1\to 0.
$$
Moreover, there is a canonical braiding $\sigma^\bullet$ on $\Omega^\bullet(A)$ such that the latter becomes graded-braided commutative \cite{DurGauge} and it follows that
\begin{equation}\label{intr:HG}
	\Omega^\bullet(B)=\Omega^\bullet(A)^{\mathrm{co}\Omega^\bullet(H)}\subseteq\Omega^\bullet(A)
\end{equation}
constitutes a graded Hopf--Galois extension \cite{DurGalois}. The \DJ ur\dj evi\'c approach to quantum principal bundles has recently been revisited in \cite{AntEmaTho}. In this proceeding we review this account, with a particular focus on new examples. We further explore the notion of quantum gauge transformations, as developed in \cite{Brzezinksi}. As in the classical setup, there is a group isomorphism relating gauge transformations and vertical automorphisms of quantum principal bundles. In order to obtain an action of gauge transformations on connections we extend the definition of quantum gauge transformations to differential forms. This turns out to be a very natural generalization of \cite{Brzezinksi}, as the extended coaction $\Delta^\bullet_A$ structures $\Omega^\bullet(A)$ as a graded comodule algebra and \eqref{intr:HG} becomes a graded Hopf--Galois extension. Such a generalization was first considered in \cite{Amilcar}. We spell out the example of the noncommutative $2$-tours and exemplify how gauge transformations act on its connections. An alternative approach to quantum gauge transformations involving (quasi)triangular Hopf algebras can be found in \cite{PaGiCh}. As possible future applications it would be interesting to investigate noncommutative differential operators and jet bundles \cite{jets,jets2}, as well as pseudo-Riemannian metrics and Levi-Civita connections \cite{BeggsMajid,PaoloThomas} in this quantum principal bundle formalism.

\subsection*{Conventions and notation}

Throughout this paper we fix a field $\Bbbk$ and denote the tensor product of $\Bbbk$-vector spaces by $\otimes$. All modules, comodules, algebras and coalgebras are understood over $\Bbbk$. If not stated otherwise, we only consider associative unital algebras and coassociative counital coalgebras. Given a coalgebra $\mathcal{C}$ with comultiplication $\Delta\colon\mathcal{C}\to\mathcal{C}\otimes\mathcal{C}$ and counit $\varepsilon\colon\mathcal{C}\to\Bbbk$, we employ Sweedler's notation for the coproduct $\Delta(c)=:c_1\otimes c_2$, where $c\in\mathcal{C}$, and, inductively 
\begin{align*}
	(\Delta\otimes\mathrm{id})(\Delta(c))
	\overset{(*)}{=}(\mathrm{id}\otimes\Delta)(\Delta(c))&=:c_1\otimes c_2\otimes c_3,\\
	(\Delta\otimes\mathrm{id}\otimes\mathrm{id})((\Delta\otimes\mathrm{id})(\Delta(c)))&=:c_1\otimes c_2\otimes c_3\otimes c_4,\\
	&\qquad\qquad\vdots
\end{align*}
where $(*)$ holds by coassociativity. Moreover, the counitality axiom reads $\varepsilon(c_1)c_2=c=c_1\varepsilon(c_2)$ for all $c\in\mathcal{C}$. The antipode of a Hopf algebra $H$ is denoted by $S\colon H\to H$ and we always assume that $S$ is invertible as a $\Bbbk$-linear map. Given an algebra $A$ the category of left $A$-modules is denoted by ${}_A\mathcal{M}$. The corresponding morphism are called left $A$-linear maps. Similarly for the category $\mathcal{M}_A$ of right $A$-modules and the category ${}_A\mathcal{M}_A$ of $A$-bimodules. Given a coalgebra $(\mathcal{C},\Delta,\varepsilon)$ we denote the category of left $\mathcal{C}$-comodules by ${}^\mathcal{C}\mathcal{M}$ and call the corresponding morphisms left $\mathcal{C}$-colinear. For a left $\mathcal{C}$-comodule $M$ with left $\mathcal{C}$-coaction $\lambda_M\colon M\to\mathcal{C}\otimes M$ we employ Sweedler's notation
\begin{align*}
	\lambda_M(m)&=:m_{-1}\otimes m_0,\\
	(\Delta\otimes\mathrm{id})(\lambda_M(m))\overset{(**)}{=}
	(\mathrm{id}\otimes\lambda_M)(\lambda_M(m))&=:m_{-2}\otimes m_{-1}\otimes m_0,\\
	&\qquad\qquad\vdots
\end{align*}
for $m\in M$, where $(**)$ holds by the coaction axiom of $\lambda_M$.
Similarly for the category $\mathcal{M}^\mathcal{C}$ of right $\mathcal{C}$-comodules and the category ${}^\mathcal{C}\mathcal{M}^\mathcal{C}$ of $\mathcal{C}$-bicomodules The short notation associated to a right $\mathcal{C}$-comodule $M$ with right $\mathcal{C}$-coaction $\Delta_M\colon M\to M\otimes\mathcal{C}$ is $\Delta_M(m)=:m_0\otimes m_1$, etc. For more information on comodules and Hopf algebras we refer to the textbook \cite{Montgomery}.

\section{Noncommutative calculi on quantum principal bundles}

In this preliminary section we recall the algebraic concept of Hopf--Galois extension and its geometric interpretation as quantum principal bundle. We then continue by discussing differential calculi on algebras, with particular focus on covariant calculi on Hopf algebras. As fundamental construction, we recall the maximal prolongation of a first order differential calculus and its explicit description in terms of invariant forms in case of covariant calculi on Hopf algebras.

\subsection{Quantum principal bundles}

Fix a Hopf algebra $(H,\Delta,\varepsilon,S)$. Recall that a \textit{right $H$-comodule algebra} is an algebra $A$, together with a right $H$-coaction $\Delta_A\colon A\to A\otimes H$, such that $\Delta_A$ is an algebra morphism. In particular,
\begin{equation}
	B:=A^{\mathrm{co}H}:=\{a\in A~|~\Delta_A(a)=a\otimes 1\}\subseteq A
\end{equation}
is a subalgebra. We call $B$ the \textit{subalgebra of coinvariant elements}. A morphism of right $H$-comodule algebras is an algebra morphism that is right $H$-colinear. Given a right $H$-comodule algebra with coinvariant subalgebra $B:=A^{\mathrm{co}H}$ we call $B\subseteq A$ a \textit{Hopf--Galois extension} \cite{KreimerTakeuchi} if the canonical map
\begin{equation}\label{GaloisMap}
	\chi\colon A\otimes_BA\to A\otimes H,\qquad\qquad
	\chi(a\otimes_Ba'):=a\Delta_A(a')
\end{equation}
is a bijection. For the restricted inverse $\tau\colon H\to A\otimes_BA$, $\tau(h):=\chi^{-1}(1\otimes h)$, called \textit{translation map}, we use the common short notation
\begin{equation}\label{translation}
	\tau(h)=:h^{\langle 1\rangle}\otimes_Bh^{\langle 2\rangle},
\end{equation}
where $h\in H$. It is well-known, see e.g. \cite{Brzezinksi}, that the equalities
\begin{align*}
	h^{\langle 1\rangle}(h^{\langle 2\rangle})_0\otimes (h^{\langle 2\rangle})_1
	&=1\otimes h,\\
	a_0(a_1)^{\langle 1\rangle}\otimes_B(a_1)^{\langle 2\rangle}
	&=1\otimes_Ba,\\
	\tau(hg)
	&=g^{\langle 1\rangle}h^{\langle 1\rangle}\otimes_Bh^{\langle 2\rangle}g^{\langle 2\rangle},\\
	h^{\langle 1\rangle}h^{\langle 2\rangle}
	&=\varepsilon(h)1,\\
	h^{\langle 1\rangle}\otimes_B(h^{\langle 2\rangle})_0\otimes(h^{\langle 2\rangle})_1
	&=(h_1)^{\langle 1\rangle}\otimes_B(h_1)^{\langle 2\rangle}\otimes h_2,\\
	(h^{\langle 1\rangle})_0\otimes_Bh^{\langle 2\rangle}\otimes(h^{\langle 1\rangle})_1
	&=(h_2)^{\langle 1\rangle}\otimes_B(h_2)^{\langle 2\rangle}\otimes S(h_1),\\
	b\tau(h)
	&=\tau(h)b
\end{align*}
hold for all $h,g\in H$, $a\in A$ and $b\in B$. 

It turns out that Hopf--Galois extensions $B=A^{\mathrm{co}H}\subseteq A$ allow for a geometric interpretation: we can understand the injection $B\hookrightarrow A$ as the algebraic counterpart of a bundle projection $P\to M$, where $A$ is the \textit{total space algebra} and $B$ the \textit{base space algebra}. The coaction of the Hopf algebra  $H$ replaces a Lie group action, thus we interpret $H$ as the \textit{structure Hopf algebra} and the \textit{quantum $H$-bundle} $B\hookrightarrow A$ is \textit{principal} by the invertibility of the Galois map \eqref{GaloisMap}. This is the point of view taken by Schneider \cite{Schneider} and in the work of \DJ ur\dj evi\'c \cite{DurII} such an algebraic structure is referred to as \textit{quantum principal bundle}. Note that there is an implicit faithfully flatness assumption, which we are discussing below, see also \cite{AntEmaTho}.
\begin{definition}[Quantum principal bundle]\label{def:QPB}
A Hopf--Galois extension $B:=A^{\mathrm{co}H}\subseteq A$ is called a \textit{quantum principal bundle} if $A$ is a faithfully flat right $B$-module.
\end{definition}
Recall that $A$ is a \textit{faithfully flat right $B$-module} if the functor $A\otimes_B\cdot\colon{}_B\mathcal{M}\to{}_\Bbbk\mathrm{Vec}$ is exact and reflects exactness. By the seminal paper \cite{Schneider}, right $B$-exactness is equivalent to left exactness (recall that we assume invertibility of the antipode and that we are working over a field $\Bbbk$), which is equivalent to the equivalence ${}_B\mathcal{M}\cong{}_A\mathcal{M}^H$ of categories. We would like to remark that in \cite{BrzMaj} quantum principal bundle is used for an $H$-comodule algebra $A$, together with first order differential calculi such that the noncommutative Atiyah sequence is exact. We comment on this in Section~\ref{sec:BrzMaj} and clarify how the two notions of quantum principal bundle relate. In the following we use Definition \ref{def:QPB} if we refer to a quantum principal bundle (QPB). 

The most immediate examples of quantum principal bundles (QPBs) are discussed below.
\begin{example}\label{ex:HopfGalois}
We fix a Hopf algebra $H$ with invertible antipode $S$.
\begin{enumerate}
\item[i.)] If we view $H$ as a right $H$-comodule algebra via the comultiplication $\Delta\colon H\to H\otimes H$, then the invariants under this coaction $H^{\mathrm{co}H}\cong\Bbbk$ are precisely the constant multiples of the unit and $\Bbbk\subseteq H$ is a Hopf--Galois extension because the antipode provides an inverse of the Galois map. Faithful flatness of $H$ as a $\Bbbk$-module is automatic, since $H$ is a $\Bbbk$-vector space. Thus, $\Bbbk\subseteq H$ is a QPB.

\item[ii.)] Consider a Lie group $G$ and a classical $G$-principal bundle $P\to M$. Then $H:=\mathscr{C}^\infty(G)$ is a Hopf algebra (where one considers the completed tensor product in the Fr\'echet topology) and by defining $B:=\mathscr{C}^\infty(M)$ and $A:=\mathscr{C}^\infty(M)$ we obtain a QPB $B\cong A^{\mathrm{co}H}\subseteq A$. Invertibility of the Glaois map \eqref{GaloisMap} follows from invertibility of the map $P\times G\to P\times_MP$. The latter is a consequence of freeness and transitivity of the $G$-action. Faithful flatness is obtained by properness of the action. This is an a posteriori motivation to consider QPBs as generalizations of classical $G$-principal bundles. More details can be found in \cite[Example~2.13]{AscBiePagSch}.

\item[iii.)] Let $A$ be a right $H$-comodule algebra and $B:=A^{\mathrm{co}H}\subseteq A$ a \textbf{cleft extension}, i.e. we assume that there is a right $H$-colinear convolution invertible map $j\colon H\to A$. Then $B\subseteq A$ is a Hopf--Galois extension, where the inverse of \eqref{GaloisMap} is given by
$$
\chi^{-1}\colon A\otimes H\to A\otimes_BA,\qquad\qquad
\chi^{-1}(a\otimes h):=aj^{-1}(h_1)\otimes_Bj(h_2).
$$
Cleft extensions are automatically faithfully flat, as shown in \cite[Part VII, Section 6]{HajacNotes}. An instance of a cleft extension is given by the noncommutative $2$-tours with $U(1)$-fibration, see \cite[Example 3.7 $ii.)$]{AntEmaTho} for more details.
\end{enumerate}
\end{example}
The important example of the quantum Hopf fibration will be discussed in Section \ref{sec:examples}.

\subsection{Differential calculi}

For convenience of the reader we briefly review the theory of noncommutative differential calculi on algebras and comodule algebras, following \cite{Woronowicz}, see also \cite{BeggsMajid,DurI,Schauenburg}.

A \textit{differential calculus} (abbreviated by DC) on an algebra $A$ is a differential graded algebra (abbreviated by DGA) $\Omega^\bullet(A)=\bigoplus\nolimits_{n\geq 0}\Omega^n(A)$, such that $\Omega^0(A)=A$ and for all $n>0$ we have 
$\Omega^n(A)=\mathrm{span}_\Bbbk\{a^0\mathrm{d}(a^1)\wedge\ldots\wedge\mathrm{d}(a^n)~|~a^0,a^1,\ldots,a^n\in A\}$. More explicitly, we have graded maps $\wedge\colon\Omega^\bullet(A)\otimes\Omega^\bullet(A)\to\Omega^\bullet(A)$ (the \textit{wedge product}) of degree $0$ and $\mathrm{d}\colon\Omega^\bullet(A)\to\Omega^\bullet(A)$ (the \textit{differential}) of degree $1$, such that $\wedge$ is associative, $\mathrm{d}^2=0$ is nilpotent and the \textit{Leibniz rule}
$$
\mathrm{d}(\omega\wedge\eta)=\mathrm{d}\omega\wedge\eta+(-1)^{|\omega|}\omega\wedge\mathrm{d}\eta
$$
holds for all $\omega,\eta\in\Omega^\bullet(A)$, where $|\omega|\in\mathbb{N}$ denotes the degree of $\omega$.

Given a Hopf algebra $A$ and a right $H$-comodule algebra $A$, we call a DC $\Omega^\bullet(A)$ on $A$ \textit{right $H$-covariant}, if $\Omega^n(A)\in{}_A\mathcal{M}_A^H$ are right $H$-covariant $A$-bimodules (i.e. $A$-bimodules and right $H$-comodules such that the coactions $\Delta_{\Omega^n(A)}\colon\Omega^n(A)\to\Omega^n(A)\otimes H$ are $A$-bilinear) and $\wedge$, $\mathrm{d}$ are right $H$-colinear maps. The latter reads
$$
\Delta_{\Omega^{|\omega|+|\eta|}(A)}(\omega\wedge\eta)=\Delta_{\Omega^{|\omega|}(A)}(\omega)\Delta_{\Omega^{|\eta|}(A)}(\eta),\qquad
\Delta_{\Omega^{|\omega|+1}}\mathrm{d}(\omega)=(\mathrm{d}\otimes\mathrm{id})(\Delta_{\Omega^{|\omega|}}(\omega))
$$
for all $\omega,\eta\in\Omega^\bullet(A)$. We denote the collection of right $H$-coactions $\Delta_{\Omega^n(A)}$ by $\Delta_{\Omega^\bullet(A)}:=\bigoplus\nolimits_{n\geq 0}\Delta_{\Omega^n(A)}$ and it structures $\Omega^\bullet(A)$ as a graded right $H$-comodule. Similarly, left $H$-covariant differential calculi and $H$-bicovariant differential calculi are defined on left $H$-comodule algebras and $H$-bicomodule algebras, respectively.

The truncation $\Omega^{\leq 1}(A)$ of a DC $\Omega^\bullet(A)$ on $A$ is called a \textit{first order differential calculus} (abbreviated by FODC) on $A$. Explicitly, a FODC on $A$ is a tuple $(\Omega^1(A),\mathrm{d})$, where $\Omega^1(A)$ is an $A$-bimodule and $\mathrm{d}\colon A\to\Omega^1(A)$ is $\Bbbk$-linear such that for all $a,b\in A$ the Leibniz rule
$$
\mathrm{d}(ab)=\mathrm{d}(a)b+a\mathrm{d}(b)
$$
holds and $\Omega^1(A)=A\mathrm{d}A$ coincides with the left $A$-module generated by the image $\mathrm{d}A$ of $\mathrm{d}$. For a right $H$-covariant FODC $\Omega^1(A)$ is assumed to be right $H$-covariant and $\mathrm{d}$ to be right $H$-colinear, in addition. Conversely, we are able to extend any FODC $(\Omega^1(A),\mathrm{d})$ on $A$ to a DC $\Omega^\bullet(A)$. This can be done trivially, i.e. by setting $\Omega^n(A)=0$ for $n>1$, or in a maximal way, using the \textit{maximal prolongation}, see e.g. \cite[Appendix B]{DurI} and \cite[Section 4]{Schauenburg}. 
\begin{proposition}
Let $A$ be an algebra and $(\Omega^1(A),\mathrm{d})$ a FODC on $A$. Then there exists a DC $\Omega^\bullet_\mathrm{max}(A)$ on $A$ such that its truncation coincides with $(\Omega^1(A),\mathrm{d})$ and such that the following universal property holds: for all differential calculi $\tilde{\Omega}^\bullet(A)$ which extend $(\Omega^1(A),\mathrm{d})$ there is a surjective morphism $\Omega^\bullet_\mathrm{max}(A)\twoheadrightarrow\tilde{\Omega}^\bullet(A)$ of differential graded algebras.

If $A$ is a right $H$-comodule algebra and $(\Omega^1(A),\mathrm{d})$ a right $H$-covariant FODC, then $\Omega^\bullet_\mathrm{max}(A)$ is right $H$-covariant.
\end{proposition}
We call $\Omega^\bullet_\mathrm{max}(A)$ the \textit{maximal prolongation} of a FODC $(\Omega^1(A),\mathrm{d})$ on $A$. Its universal property says that all differential calculi extending $(\Omega^1(A),\mathrm{d})$ are quotient differential calculi of $\Omega^\bullet_\mathrm{max}(A)$ and thus $\Omega^\bullet_\mathrm{max}(A)$ is the maximal extension. There is an explicit description of $\Omega^\bullet_\mathrm{max}(A)$, which we sketch in the following: given a FODC $(\Omega^1(A),\mathrm{d})$ on $A$, consider the tensor algebra $T\Omega^1(A):=A\oplus\Omega^1(A)\oplus(\Omega^1(A)\otimes_A\Omega^1(A))\oplus\ldots$ of $\Omega^1(A)$ and quotient it by the ideal generated by elements
$$
\mathrm{d}a^i\otimes_A\mathrm{d}b^i,\qquad\qquad
\text{where }a^i,b^i\in A,\text{ such that }a^i\mathrm{d}b^i=0.
$$
By construction, the differential extends to the quotient. In the following we usually suppress the subscript of $\Omega^\bullet_\mathrm{max}(A)$.

\subsection{Covariant calculi on Hopf algebras}

We recall a special case of the previous section, namely the description of covariant and bicovariant differential calculi on Hopf algebras. For this we follow again the original source \cite{Woronowicz}, as well as \cite{BeggsMajid,DurI,Schauenburg}.

Let $H$ be a Hopf algebra with invertible antipode. Given a left $H$-covariant (in the following \textit{left covariant} for short) FODC $(\Omega^1(H),\mathrm{d})$ on $H$, we denote the left $H$-coaction on $\Omega^1(H)$ by $\lambda_{\Omega^1(H)}\colon\Omega^1(H)\to H\otimes\Omega^1(H)$ and the vector subspace of invariant elements by
$$
\Lambda^1:={}^{\mathrm{co}H}\Omega^1(H):=\{\omega\in\Omega^1(H)~|~\lambda_{\Omega^1(H)}(\omega)=1\otimes\omega\}.
$$
There is a surjection from the kernel of the counit $H^+:=\ker\varepsilon\subset H$ to $\Lambda^1$, given by
\begin{equation}
\varpi\colon H^+\twoheadrightarrow\Lambda^1,\qquad\qquad
\varpi(h):=S(h_1)\mathrm{d}(h_2),
\end{equation}
called the \textit{Cartan--Maurer form}. If we endow $\Lambda^1$ with the adjoint right $H$-action $\vartheta\cdot h:=S(h_1)\vartheta h_2$, then $\varpi$ becomes right $H$-linear. Denoting the kernel of $\varpi$ by $I:=\ker\varpi\subseteq H^+$, there are bijections
\begin{equation}
	\Omega^1(H)\cong H\otimes\Lambda^1\cong H\otimes H^+/I
\end{equation}
of left $H$-covariant $H$-bimodules. If we endow $H\otimes H^+/I$ with the differential $\mathrm{d}'\colon H\to H\otimes H^+/I$, $\mathrm{d}'(h):=(\mathrm{id}\otimes\pi)(\Delta(h)-h\otimes 1)$, where $\pi\colon H^+\to H^+/I$ is the projection, the above becomes an isomorphism of left covariant FODCi. This leads to the following classification theorem of Woronowicz \cite{Woronowicz}.
\begin{theorem}\label{thm:Wor}
For a Hopf algebra $H$ there is a one-to-one correspondence between left covariant FODCi on $H$ and right $H$-ideals of $H^+$. A left covariant FODC on $H$ is bicovariant if and only if the corresponding ideal $I\subseteq H^+$ is closed under the adjoint coaction, i.e. $\mathrm{Ad}(I)\subseteq I\otimes H$, where $\mathrm{Ad}(h):=h_2\otimes S(h_1)h_3$ for all $h\in H$.
\end{theorem}
Given a left covariant FODC $(\Omega^1(H),\mathrm{d})$ on $H$, we are able to describe its maximal prolongation $\Omega^\bullet(H)$ in terms of the left coinvariant $1$-forms $\Lambda^1\cong H^+/I$. As described in \cite[Proposition 2.31]{BeggsMajid}, the fundamental theorem of Hopf modules gives an isomorphism $\Omega^\bullet(H)\cong H\otimes\Lambda^\bullet$, where $\Lambda^\bullet:={}^{\mathrm{co}H}\Omega^\bullet(H)$ is isomorphic to the free algebra generated by $\Lambda^1$, modulo the relations
$$
\wedge\circ(\varpi\circ\pi_\varepsilon\otimes\varpi\circ\pi_\varepsilon)\Delta(I)=0,
$$
where $\pi_\varepsilon\colon H\to H^+$ is the canonical projection. In particular, the \textit{Cartan--Maurer equation}
\begin{equation}
	\mathrm{d}\varpi(\pi_\varepsilon(h))
	+\varpi(\pi_\varepsilon(h_1))\wedge\varpi(\pi_\varepsilon(h_2))
	=0
\end{equation}
holds for all $h\in H$.

Explicit examples of covariant calculi will be discussed in Section \ref{sec:examples}. There, we are going to encounter the $3$-dimensional covariant calculus on the quantum Hopf fibration, as well as covariant calculi on smash product algebras and a calculus on the noncommutative $2$-torus.

\section{Complete differential calculi}

In this section we equip quantum principal bundles (QPBs) $B:=A^{\mathrm{co}H}\subseteq A$ with differential structures. This will be done in a way that allows many algebro-geometric features to emerge, such as vertical, horizontal and base forms. Moreover, the Hopf--Galois structure will automatically extend to the level of differential forms and there is a canonical braiding on the level of total space forms such that the latter become graded-braided commutative. We elaborate on the papers \cite{DurII} and \cite{BrzMaj}, following the recent account \cite{AntEmaTho}.

\subsection{Total space forms, basic forms, vertical forms and horizontal forms}

Throughout this section we fix a Hopf algebra $H$ with invertible antipode $S$ and a bicovariant FODC $\Omega^1(H)\cong H\otimes\Lambda\cong H\otimes H^+/I$ on $H$. Its maximal prolongation is denoted by $\Omega^\bullet(H)$.
\begin{lemma}\label{lem:gradedHopf}
The bicovariant differential calculus $\Omega^\bullet(H)$ is a differential graded Hopf algebra with graded comultiplication, counit and antipode determined by
\begin{align*}
\Delta^\bullet\colon\Omega^\bullet(H)\to\Omega^\bullet(H)\otimes\Omega^\bullet(H),\qquad
&\Delta^\bullet(h^0\mathrm{d}h^1\wedge\ldots\mathrm{d}h^n)=\Delta(h^0)\mathrm{d}_\otimes\Delta(h^1)\ldots\mathrm{d}_\otimes\Delta(h^n)\\
\varepsilon^\bullet\colon\Omega^\bullet(H)\to\Bbbk,\qquad
&\varepsilon^\bullet(h^0\mathrm{d}h^1\wedge\ldots\mathrm{d}h^n)=0,\\
S^\bullet\colon\Omega^\bullet(H)\to\Omega^\bullet(H),\qquad
&S^\bullet(h^0\mathrm{d}h^1\wedge\ldots\mathrm{d}h^n)=\mathrm{d}(S(h^n))\wedge\ldots\wedge\mathrm{d}(S(h^1))S(h^0)
\end{align*}
for all $n>0$ and $h^0,h^1,\ldots,h^n\in H$, where $\mathrm{d}_\otimes$ denotes the differential of the tensor product algebra $H\otimes H$.
\end{lemma}
For a proof of this lemma we refer to \cite[Lemma 5.4]{Schauenburg}. Note that the multiplication on $\Omega^\bullet(H)\otimes\Omega^\bullet(H)$ involves a sign, for example
\begin{align*}
	\Delta(h)\mathrm{d}_\otimes\Delta(g)\mathrm{d}_\otimes\Delta(k)
	&=(h_1\otimes h_2)(\mathrm{d}g_1\otimes g_2+g_1\otimes\mathrm{d}g_2)(\mathrm{d}k_1\otimes k_2+k_1\otimes\mathrm{d}k_2)\\
	&=h_1\mathrm{d}g_1\wedge\mathrm{d}k_1\otimes h_2g_2k_2
	+h_1\mathrm{d}(g_1)k_1\otimes h_2g_2\mathrm{d}k_2\\
	&\quad-h_1g_1\mathrm{d}k_1\otimes h_2\mathrm{d}(g_2)k_2
	+h_1g_1k_1\otimes h_2\mathrm{d}g_2\wedge\mathrm{d}k_2
\end{align*}
for all $h,g,k\in H$. 

To rephrase Lemma \ref{lem:gradedHopf} in other words: for every bicovariant FODC the comultiplication $\Delta\colon H\to H\otimes H$ (as well as the counit and antipode) extends to a morphism of differential graded algebras (DGAs) $\Delta^\bullet\colon\Omega^\bullet(H)\to\Omega^\bullet(H)\otimes\Omega^\bullet(H)$ on the maximal prolongation $\Omega^\bullet(H)$. In this way, the Hopf algebra structure of $H$ amplifies to the level of differential calculus. Keeping this in mind, the following definition, proposed by \DJ ur\dj evi\'c in \cite{DurII}, is natural.
\begin{definition}\label{def:complete}
Let $\Omega^\bullet(H)$ be the maximal prolongation of a bicovariant FODC on $H$ and $B:=A^{\mathrm{co}H}\subseteq A$ a quantum principal bundle. We call a DC $\Omega^\bullet(A)$ on $A$ \textbf{complete}, if the right $H$-coaction $\Delta_A\colon A\to A\otimes H$ extends to a morphism  
\begin{equation}
	\Delta_A^\bullet\colon\Omega^\bullet(A)\to\Omega^\bullet(A)\otimes\Omega^\bullet(H)
\end{equation}
of DGAs.
\end{definition}
In particular, the bicovariant DC $\Omega^\bullet(H)$ is complete for the QPB $\Bbbk\cong H^{\mathrm{co}H}\subseteq H$.

In the following we use the short notations
\begin{align}
	\Delta^\bullet(\omega)&=:\omega_{[1]}\otimes\omega_{[2]}\in\Omega^\bullet(H)\otimes\Omega^\bullet(H),\\
	\Delta_A^\bullet(\eta)&=:\eta_{[0]}\otimes\eta_{[1]}\in\Omega^\bullet(A)\otimes\Omega^\bullet(H)
\end{align}
for the coproduct of $\omega\in\Omega^\bullet(H)$ and the extended coaction of $\eta\in\Omega^\bullet(A)$ for a complete differential calculus $\Omega^\bullet(A)$, respectively.

\begin{definition}
Consider a complete differential calculus $\Omega^\bullet(A)$ on a QPB $B:=A^{\mathrm{co}H}\subseteq A$. We call
\begin{enumerate}
\item[i.)] $\Omega^\bullet(A)$ the \textbf{total space forms},

\item[ii.)] $\Omega^\bullet(B):=\Omega^\bullet(A)^{\mathrm{co}\Omega^\bullet(H)}=\{\omega\in\Omega^\bullet(A)~|~\Delta_A^\bullet(\omega)=\omega\otimes 1\}\subseteq\Omega^\bullet(A)$ the \textbf{basic forms},

\item[iii.)] $\mathrm{ver}^\bullet:=A\otimes{}^{\mathrm{co}H}\Omega^\bullet(H)=A\otimes\Lambda^\bullet$ the \textbf{vertical forms},

\item[iv.)] $\mathrm{hor}^\bullet:=(\Delta_A^\bullet)^{-1}(\Omega^\bullet(A)\otimes H)\subseteq\Omega^\bullet(A)$ the \textbf{horizontal forms}
\end{enumerate}
of the QPB.
\end{definition}
The properties and structures of the previously defined forms are summarized in the following proposition. For proofs of the statements we refer to \cite{AntEmaTho}.
\begin{proposition}
For a complete calculus $\Omega^\bullet(A)$ it follows that
\begin{enumerate}
\item[i.)] $\Omega^\bullet(A)$ is a (graded) right $\Omega^\bullet(H)$-comodule algebra with (graded) right $\Omega^\bullet(H)$-coaction $\Delta_A^\bullet\colon\Omega^\bullet(A)\to\Omega^\bullet(A)\otimes\Omega^\bullet(H)$. In particular, the DC $\Omega^\bullet(A)$ is right $H$-covariant.

\item[ii.)] $\Omega^\bullet(B)\subseteq\Omega^\bullet(A)$ is a differential graded subalgebra.

\item[iii.)] $\mathrm{ver}^\bullet=A\otimes\Lambda^\bullet$ is a differential calculus on $A$ with respect to the wedge product and differential
\begin{equation}
\begin{split}
(a\otimes\vartheta)\wedge(a'\otimes\vartheta')
&:=aa'_0\otimes S(a'_1)\vartheta a'_2\wedge\vartheta',\\
\mathrm{d}(a\otimes\vartheta)
&:=a\otimes\mathrm{d}\vartheta+a_0\otimes\varpi(\pi_\varepsilon(a_1))\wedge\vartheta,
\end{split}
\end{equation}
where $a,a'\in A$ and $\vartheta,\vartheta'\in\Lambda^\bullet$. It is a complete calculus with differential graded extension $\Delta_v^\bullet\colon\mathrm{ver}^\bullet\to\mathrm{ver}^\bullet\otimes\Omega^\bullet\colon\mathrm{ver}^\bullet\to\mathrm{ver}^\bullet\otimes\Omega^\bullet(H)$ determined by the commutative diagram
\begin{equation}
\begin{tikzcd}
\Omega^\bullet(A) \arrow{rr}{\pi_v} \arrow{d}[swap]{\Delta_A^\bullet} 
& & \mathrm{ver}^\bullet \arrow{d}{\Delta_v^\bullet}\\
\Omega^\bullet(A)\otimes\Omega^\bullet(H) \arrow{rr}{\pi_v\otimes\mathrm{id}}
& & \mathrm{ver}^\bullet\otimes\Omega^\bullet(H)
\end{tikzcd}
\end{equation}
where 
\begin{equation}
	\pi_v(a^0\mathrm{d}a^1\wedge\ldots\wedge\mathrm{d}a^n):=a_0^0a_0^1\ldots a_0^n\otimes S(a_1^0a_1^1\ldots a_1^n)a_2^0\mathrm{d}a_2^1\wedge\ldots\mathrm{d}a_2^n
\end{equation}
for all $a^0,a^1,\ldots,a^n\in A$.

\item[iv.)] $\mathrm{hor}^\bullet\subseteq\Omega^\bullet(A)$ is a graded right $H$-comodule subalgebra.
\end{enumerate}
\end{proposition}
Not that $\Omega^\bullet(B)$ is \textit{not} a DC in general. In fact, there are tautological examples where $\Omega^1(B)$ is not generated in degree $0$. However, in all explicit examples we encounter in this article the base forms $\Omega^\bullet(B)$ are isomorphic to the pullback calculus $\Omega^\bullet(B)\subseteq\Omega^\bullet(A)$ induced by the algebra inclusion $B\hookrightarrow A$.

\subsection{First order completeness and the Brzezi\'nski--Majid approach}\label{sec:BrzMaj}

In first order, horizontal forms inject in total space forms and the latter surject on vertical forms. The resulting exact sequence is the noncommutative analogue of the Atiyah sequence.
\begin{proposition}\label{prop:Atiyah}
Given a compete calculus $\Omega^\bullet(A)$ on a QPB $B:=A^{\mathrm{co}H}\subseteq A$, the noncommutative Atiyah sequence
\begin{equation}\label{Atiyah}
0\to\mathrm{hor}^1\hookrightarrow\Omega^1(A)\overset{\pi_v}{\twoheadrightarrow}\mathrm{ver}^1\to 0
\end{equation}
is exact in the category ${}_A\mathcal{M}_A^H$ of right $H$-covariant $A$-bimodules.
\end{proposition}
Obviously, it is sufficient to require completeness up to first order, rather than completeness for all degrees as in Definition \ref{def:complete}, in order to obtain the exact sequence \eqref{Atiyah}. We call a differential calculus $\Omega^\bullet(A)$ on a QPB $B:=A^{\mathrm{co}H}\subseteq A$ \textit{first order complete}, if $\Delta_A\colon A\to A\otimes H$ is $1$-differential, i.e. if $\Omega^1(A)$ is right $H$-covariant and $\pi_v\colon\Omega^1(A)\to A\otimes\Lambda^1$, $\pi_v(a\mathrm{d}a'):=aa'_0\otimes S(a'_1)\mathrm{d}a'_2$ is well-defined. In this case, 
$$
\Delta_A^1=\Delta_{\Omega^1(A)}+\mathrm{ver}^{0,1}\colon\Omega^1(A)\to(\Omega^1(A)\otimes H)\oplus(A\otimes\Omega^1(H))
$$ 
is the differential of $\Delta_A$, where $\mathrm{ver}^{0,1}\colon\Omega^1(A)\to A\otimes\Omega^1(H)$, $a\mathrm{d}a'\mapsto a_0a'_0\otimes a_1\mathrm{d}_H(a'_1)$.
\begin{remark}
In the seminal paper \cite{BrzMaj}, Brzezi\'nski and Majid define quantum principal bundles on right $H$-comodule algebras $A$ as
a covariant FODC $\Omega^1(H)$ on $H$ and a right $H$-covariant FODC $\Omega^1(A)$ on $A$, such that
\begin{equation}
	0\to A\Omega^1(B)A\hookrightarrow\Omega^1(A)\overset{\pi_v}{\twoheadrightarrow}\mathrm{ver}^1\to 0
\end{equation}
is exact, where $\Omega^1(B)\subseteq\Omega^1(A)$ is the pullback calculus $B\hookrightarrow A$. In case $B:=A^{\mathrm{co}H}\subseteq A$ is a faithfully flat Hopf--Galois extension, i.e. a QPB according to Definition \ref{def:QPB}, then any quantum principal bundle in the sense of Brzezi\'nski--Majid is first order complete. Conversely, a first order complete calculus is a quantum principal bundle in the sense of Brzezi\'nski--Majid if and only if $A\Omega^1(B)A=\mathrm{hor}^1$.
\end{remark}
So according to the previous remark, quantum principal bundles in the sense of Brzezi\'nski--Majid are effectively equivalent to first order complete calculi. In the next proposition we show that any first order complete calculus extends to a complete calculus on the maximal prolongation. This shows that the \DJ ur\dj evi\'c approach to complete calculi is effectively an extension of the Brzezi\'nski--Majid approach to quantum principal bundles.
\begin{proposition}\label{prop:maxprolcomp}
Let $B:=A^{\mathrm{co}H}\subseteq A$ be a QPB, i.e. a faithfully flat Hopf--Galois extension. If $(\Omega^1(A),\Omega^1(H))$ is a first order complete calculus on $A$, then the maximal prolongations $(\Omega^\bullet(A),\Omega^\bullet(H))$ form a complete calculus on $A$.
\end{proposition}
For a proof we refer to \cite[Proposition 3.31]{AntEmaTho}.

\subsection{Graded Hopf--Galois extension and braidings}

We recall that Hopf--Galois extensions are endowed with a canonical braiding, the \DJ ur\dj evi\'c braiding, such that the comodule algebra becomes braided-commutative. In continuation, we show that the Hopf--Galois property and the braiding amplify to the level of differential forms of complete calculi. We follow \cite{DurGauge,DurGalois}, see also \cite[Section 5]{AntEmaTho}.

Given a Hopf--Galois extension $B:=A^{\mathrm{co}H}\subseteq A$ we can use the vector space isomorphism $\chi\colon A\otimes_BA\xrightarrow{\cong}A\otimes H$ and pull back the tensor product multiplication on $A\otimes H$ in order to obtain an associative unital product
\begin{equation}
\begin{split}
	(A\otimes_BA)\otimes_B(A\otimes_BA)&\to A\otimes_BA\\
	(a\otimes_Ba')\otimes_B(c\otimes_Bc')&\mapsto a\sigma(a'\otimes_Bc)c'
\end{split}
\end{equation}
on $A\otimes_BA$, where $\sigma\colon A\otimes_BA\to A\otimes_BA$ denotes the \textit{\DJ ur\dj evi\'c braiding}. The latter is the $B$-bilinear isomorphism defined by
\begin{equation}
	\sigma(a\otimes c):=a_0c\tau(a_1)=a_0c(a_1)^{\langle 1\rangle}\otimes_B(a_1)^{\langle 2\rangle}
\end{equation}
for all $a,c\in A$. It satisfies the braid equation $\sigma_{12}\sigma_{23}\sigma_{12}=\sigma_{23}\sigma_{12}\sigma_{23}$ in $A^{\otimes_B3}$ and
\begin{align*}
	\sigma(m_A\otimes_B\mathrm{id}_A)&=(\mathrm{id}_A\otimes_Bm_A)\sigma_{12}\sigma_{23},\\
	\sigma(\mathrm{id}_A\otimes_Bm_A)&=(m_A\otimes_B\mathrm{id}_A)\sigma_{23}\sigma_{12},
\end{align*}
where $m_A\colon A\otimes_BA\to A$ denotes the multiplication. Furthermore, $A$ is braided-commutative with respect to $\sigma$, i.e.
\begin{equation}
	m_A\circ\sigma=m_A
\end{equation}
holds as an equation of endomorphisms of $A\otimes_BA$. More details and explicit examples of the \DJ ur\dj evi\'c braiding are spelled out in Section \ref{sec:examples} and \cite[Section 6]{AntEmaTho}.

For a complete calculus $\Omega^\bullet(A)$ on a QPB $B:=A^{\mathrm{co}H}\subseteq A$ we have, by assumption, the (unique) extension $\Delta_A^\bullet\colon\Omega^\bullet(A)\to\Omega^\bullet(A)\otimes\Omega^\bullet(H)$ of the right $H$-coaction on $A$ as a morphism of DGAs. In particular, there is a graded extension
\begin{equation}\label{gradedGalois}
\begin{split}
	\chi^\bullet\colon\Omega^\bullet(A)\otimes_{\Omega^\bullet(B)}\Omega^\bullet(A)&\to\Omega^\bullet(A)\otimes\Omega^\bullet(H),\\
	\omega\otimes_{\Omega^\bullet(B)}\eta&\mapsto\omega\wedge\Delta_A^\bullet(\eta)=\omega\wedge\eta_{[0]}\otimes\eta_{[1]}
\end{split}
\end{equation}
of the Galois map \eqref{GaloisMap}. It turns out that \eqref{gradedGalois} is automatically invertible, leading to the following result, shown in \cite[Proposition 2]{DurGalois} and \cite[Theorem 5.3]{AntEmaTho}, see also \cite[Satz 5.5.6]{SchauenburgThesis}.
\begin{theorem}
Let $\Omega^\bullet(A)$ be a complete calculus on a QPB $B:=A^{\mathrm{co}H}\subseteq A$. Then
\begin{equation}
	\Omega^\bullet(B)=\Omega^\bullet(A)^{\mathrm{co}\Omega^\bullet(H)}\subseteq\Omega^\bullet(A)
\end{equation}
is a graded Hopf--Galois extension, i.e. \eqref{gradedGalois} is invertible.
\end{theorem}
Similarly to the non-graded case, one defines a graded analogue
\begin{equation}
	\tau^\bullet:=(\chi^\bullet)^{-1}\colon\Omega^\bullet(H)\to\Omega^\bullet(A)\otimes_{\Omega^\bullet(B)}\Omega^\bullet(A)
\end{equation}
of the translation map, which satisfies equations similar to the ones of $\tau$ (see \cite[Proposition]{AntEmaTho} for more details). Moreover, we can use the extended translation map $\tau^\bullet(\omega)=:\omega^{[1]}\otimes_{\Omega^\bullet(B)}\omega^{[2]}$, for $\omega\in\Omega^\bullet(H)$, to define the graded extension
\begin{equation}\label{extendedBraid}
\begin{split}
	\sigma^\bullet\colon\Omega^\bullet(A)\otimes_{\Omega^\bullet(B)}\Omega^\bullet(A)&\to\Omega^\bullet(A)\otimes_{\Omega^\bullet(B)}\Omega^\bullet(A)\\
	\omega\otimes_{\Omega^\bullet(B)}\eta&\mapsto
	(-1)^{|\eta||\omega_{[1]}|}\omega_{[0]}\wedge\eta\wedge\tau^\bullet(\omega_{[1]})
\end{split}
\end{equation}
of the \DJ ur\dj evi\'c braiding. We give the properties of $\sigma^\bullet$, referring \cite[Proposition 3.1]{DurGauge} and \cite[Proposition 5.6]{AntEmaTho} for a proof.
\begin{proposition}
For every complete calculus $\Omega^\bullet(A)$ the following statements hold true.
\begin{enumerate}
\item[i.)] The extension $\sigma^\bullet$ of the \DJ ur\dj evi\'c, defined in \eqref{extendedBraid}, is a $\Omega^\bullet(B)$-bilinear isomorphism.

\item[ii.)] $\sigma^\bullet$ satisfies the braid relation $\sigma^\bullet_{12}\sigma^\bullet_{23}\sigma^\bullet_{12}=\sigma^\bullet_{23}\sigma^\bullet_{12}\sigma^\bullet_{23}$ in $\Omega^\bullet(A)^{\otimes_{\Omega^\bullet(B)}3}$.

\item[iii.)] One has
\begin{align*}
	\sigma^\bullet(\wedge\otimes_{\Omega^\bullet(B)}\mathrm{id}_{\Omega^\bullet(A)})
	&=(\mathrm{id}_{\Omega^\bullet(A)}\otimes_{\Omega^\bullet(B)}\wedge)\sigma^\bullet_{12}\sigma^\bullet_{23},\\
	\sigma^\bullet(\mathrm{id}_{\Omega^\bullet(A)}\otimes_{\Omega^\bullet(B)}\wedge)
	&=(\wedge\otimes_{\Omega^\bullet(B)}\mathrm{id}_{\Omega^\bullet(A)})\sigma^\bullet_{23}\sigma^\bullet_{12},
\end{align*}
where $\wedge\colon\Omega^\bullet(A)\otimes_{\Omega^\bullet(B)}\Omega^\bullet(A)\to\Omega^\bullet(A)$ denotes the wedge product of $\Omega^\bullet(A)$.

\item[iv.)] $\Omega^\bullet(A)$ is graded-braided-commutative with respect to $\sigma^\bullet$, i.e.
\begin{equation}\label{gradedbraidedcommutative}
	\wedge\circ\sigma^\bullet=\wedge
\end{equation}
holds as an equation of maps $\Omega^\bullet(A)\otimes_{\Omega^\bullet(B)}\Omega^\bullet(A)\to\Omega^\bullet(A)$.
\end{enumerate}
\end{proposition}
Note that $\sigma^\bullet$ does not square to the identity in general, i.e. it is not a symmetric braiding. Further note that for a commutative algebra, with a graded-commutative differential calculus on it, $\sigma^\bullet$ reduces to the graded flip isomorphism and \eqref{gradedbraidedcommutative} reflects the usual graded-commutativity of forms, where the sign is absorbed by \eqref{extendedBraid}. For explicit examples of $\sigma^\bullet$ we refer to interested reader to the following section and \cite[Section 6]{AntEmaTho}.

\subsection{Examples of complete calculi}\label{sec:examples}

In this section we discuss several (classes of) examples of complete differential calculi and their emerging geometric structures.

\paragraph{Hopf algebras:} Every Hopf algebra $H$ is a right $H$-comodule algebra over itself with coaction given by the comultiplication. The corresponding coinvariants $B:=H^{\mathrm{co}H}\cong\Bbbk$ are the constant multiples of the identity. Moreover, $\Bbbk\subseteq H$ is Hopf--Galois because in this case the antipode provides an inverse to the Galois map \eqref{GaloisMap}. Faithful flatness is automatic since $H$ is a $\Bbbk$-vector space. Thus, $\Bbbk\subseteq H$ is a QPB. Using Proposition \ref{prop:maxprolcomp} and Lemma \ref{lem:gradedHopf} one verifies the following result, taken from \cite[Section 6.1]{AntEmaTho}.
\begin{proposition}
Let $H$ be a Hopf algebra, viewed as a QPB $\Bbbk\cong H^{\mathrm{co}H}\subseteq H$. Then for any bicovariant FODC $(\Omega^1(H),\mathrm{d})$ on $H$, the maximal prolongation $\Omega^\bullet(H)$ is a complete calculus. Moreover,
\begin{enumerate}
\item[i.)] basic forms and horizontal forms are trivial, while $\mathrm{ver}^\bullet\cong\Omega^\bullet(H)$.

\item[ii.)] the extended translation map reads
\begin{equation}
	\tau^\bullet(\omega)=S^\bullet(\omega_{[1]})\otimes\omega_{[2]}
\end{equation}
for all $\omega\in\Omega^\bullet(H)$.

\item[iii.)] the extended \DJ ur\dj evi\'c braiding
\begin{equation}
	\sigma^\bullet(\omega\otimes\eta)=(-1)^{|\eta||\omega_{[2]}|}\omega_{[1]}\wedge\eta\wedge\tau^\bullet(\omega_{[2]}),
\end{equation}
where $\omega,\eta\in\Omega^\bullet(H)$, coincides with the (graded) Yetter--Drinfel'd braiding.
\end{enumerate}
\end{proposition}
From Theorem \ref{thm:Wor} we obtain a multitude of bicovariant calculi and thus many examples of complete calculi. Of particular interest are the bicovariant calculi on matrix quantum groups \cite{Jurco}.

\paragraph{The noncommutative algebraic 2-torus.} 
In this section we discuss the quantum principal bundle associated with the noncommutative algebraic 2-torus. This is given by the cleft extension $\mathcal{O}_{\theta}(\mathbb{T}^{2})^{\mathrm{co}H} \subseteq \mathcal{O}_{\theta}(\mathbb{T}^{2}) = \mathbb{C}[u,v,u^{-1},v^{-1}]/\langle vu - e^{i\theta} uv \rangle$, where $\theta \in \mathbb{R} \setminus \{0\}$ and $H = \mathcal{O}(U(1)) = \mathbb{C}[t, t^{-1}]$. The subalgebra of coinvariant elements under the right $H$-coaction is $B = \{(uv)^k \mid k \in \mathbb{Z}\}$.

We define on $A = \mathcal{O}_\theta(\mathbb{T}^2)$ a differential calculus with $\Omega^1(A) := \mathrm{span}_A\{\mathrm{d}u, \mathrm{d}v\}$ and $\Omega^2(A) := \mathrm{span}_A\{\mathrm{d}u \wedge \mathrm{d}v\}$, where the commutation relations for the calculus follow from those of $A$. Moreover, $\Omega^k(A) = 0$ for $k > 2$. We equip $H$ with its classical bicovariant FODC $\Omega^{1}(H) = \mathrm{span}_{H}\{\mathrm{d}t\}$ with trivial commutation relations alongside elements of $A$. The right $H$-coaction $\Delta_A \colon A \to A \otimes H$ lifts to 
\[
\Delta_A^1 := \Delta_{\Omega^1(A)} + \mathrm{ver}^{0,1} \colon \Omega^1(A) \to (A \otimes \Omega^1(H)) \oplus (\Omega^1(A) \otimes H),
\]
where
\[
\begin{aligned}
	\Delta_{\Omega^1(A)}(\mathrm{d}u) &= \mathrm{d}u \otimes t, & \qquad
	\Delta_{\Omega^1(A)}(\mathrm{d}v) &= \mathrm{d}v \otimes t^{-1}, \\
	\mathrm{ver}^{0,1}(\mathrm{d}u) &= u \otimes \mathrm{d}t, & \qquad
	\mathrm{ver}^{0,1}(\mathrm{d}v) &= v \otimes \mathrm{d}t^{-1}.
\end{aligned}
\]
Accordingly, the right $H$-coaction on $A$ extends to $\Omega^2(A)$.

\begin{theorem}[{\cite[Theorem 6.1]{AntEmaTho}}]
	The coaction $\Delta_A \colon A \to A \otimes H$ extends to $\Omega^{\bullet}(A)$ as a morphism of differential graded algebras. Therefore, $\Omega^{\bullet}(A)$ is a complete differential calculus. Moreover, $\Omega^\bullet(B) \subseteq \Omega^\bullet(A)$ is the pullback differential calculus.
\end{theorem}
We summarize the braiding, both, at the level of algebra and differential calculus. Since $B \subseteq A$ is a cleft extension with cleaving map $j:H\to A$, the translation map is easily understood in terms of the $j:H\to A$ and the corresponding convolution inverse, namely $\tau(h) = j^{-1}(h_1) \otimes_B j(h_2)$. The braiding $\sigma: A \otimes_B A \to A \otimes_B A$ acts accordingly at the level of algebra, and moreover we extend it naturally to differential forms. In \cite[Section 6.2]{AntEmaTho} we conclude that in this example $\sigma\colon A\otimes_BA\to A\otimes_BA$ squares to the identity. Furthermore, $\sigma^{\bullet}\colon\Omega^{\bullet}(A)\otimes_{\Omega^{\bullet}(B)}\Omega^{\bullet}(A)\rightarrow\Omega^{\bullet}(A)\otimes_{\Omega^{\bullet}(B)}\Omega^{\bullet}(A)$ squares to the identity on generators and maps generators to generators. Thus, $\sigma^{\bullet}\colon\Omega^{\bullet}(A)\otimes_{\Omega^{\bullet}(B)}\Omega^{\bullet}(A)\rightarrow\Omega^{\bullet}(A)\otimes_{\Omega^{\bullet}(B)}\Omega^{\bullet}(A)$ squares to the identity by the hexagon relations. Explicit formulas can be found in the aforementioned reference.

\paragraph{The quantum Hopf fibration:}
In this section we describe the quantum Hopf fibration as an example of quantum principal bundle with a complete differential calculus. We consider the quantum group $\mathcal{O}_q(\text{SU}(2))$ as a right $\mathcal{O}(\text{U}(1))$-comodule algebra. 
The algebra $A = \mathcal{O}_q(\text{SU}(2))$ is freely generated by the elements $\alpha, \beta, \gamma, \delta$, subject to $q$-deformed commutation relations 
\begin{equation}\nonumber 
	\begin{aligned}
		\beta \alpha  &= q \alpha \beta, \quad \gamma \alpha = q \alpha \gamma, \quad \delta \beta  = q \beta \delta,\\ 
		\quad \delta \gamma & = q \gamma \delta, \quad \gamma \beta  = \beta \gamma, \quad \delta \alpha - \alpha \delta = (q - q^{-1}) \beta \gamma
	\end{aligned}
\end{equation}
alongside the quantum determinant condition $\alpha \delta - q^{-1} \beta \gamma = 1$. The algebra $H = \mathcal{O}(\text{U}(1))$ is the Hopf algebra of rational polynomials on one variable $t$. 
There is a right $H$-action $\Delta_A: A\longrightarrow A \otimes H$, determined on generators by
\begin{equation}\nonumber 
	\begin{aligned}
		\begin{pmatrix} 
			\alpha & \beta \\ 
			\gamma & \delta 
		\end{pmatrix} &\mapsto 
		\begin{pmatrix} 
			\alpha \otimes t & \beta \otimes t^{-1} \\ 
			\gamma \otimes t & \delta \otimes t^{-1} 
		\end{pmatrix}
	\end{aligned}
\end{equation}
and extended as an algebra morphism.
The Podle\'s sphere \cite{Podles:1987wd}, denoted $B = A^{\text{co}H}$, consists of elements invariant under this coaction and it is generated by $B_+ = \alpha \beta$, $B_- = \gamma \delta$, and $B_0 = \gamma \beta$, satisfying the relations  
\begin{equation}\nonumber
	B_- B_0 = q^2 B_0 B_-, \quad B_- B_+ = q^2 B_0 (1 - q^2 B_0), \quad B_+ B_- = B_0 (1 - B_0).
\end{equation}

We construct a  differential calculus $\Omega^\bullet(A)$ on $A$ as in \cite[Example 2.32]{BeggsMajid}. Starting with the FODC $\Omega^1(A)$ defined as the free left $A$-module generated by elements 
\begin{equation}\nonumber
	e^+ = q^{-1} \alpha d\gamma - q^{-2} \gamma d\alpha,\quad e^- = \delta d\beta - q \beta d\delta, \quad e^0 = \delta d\alpha - q \beta d\gamma,    \end{equation}
with commutation rules $e^\pm f = q^{|f|} f e^\pm$ and $e^0 f = q^{2|f|} f e^0$ for $f \in \{\alpha, \beta, \gamma, \delta\}$, where $|f|$ is the degree defined by $\Delta_A(f) = f \otimes t^{|f|}$. The coaction $\Delta_{A}:A\to A\otimes H$ extends to $\Omega^1(A)$ as $\Delta_{A}^{1}:=\Delta_{\Omega^{1}(A)}+\text{ver}^{0,1}:\Omega^{1}(A)\to (\Omega^{1}(A)\otimes H) \oplus (A\otimes \Omega^{1}(H))$, with
\begin{equation} \nonumber
	\Delta_{\Omega^1}(e^\pm) = e^\pm \otimes t^{\pm 2}, \quad \Delta_{\Omega^1}(e^0) = e^0 \otimes 1,\quad \text{ver}^{0,1}(e^0) = 1 \otimes t^{-1} dt,\quad \text{ver}^{0,1}(e^\pm) = 0. 
\end{equation}
Higher-order calculi are defined recursively: $\Omega^2(A)$ is spanned by $e^\pm \wedge e^0$ and $e^+ \wedge e^-$, with relations,  while $\Omega^3(A)$ is generated by $e^+ \wedge e^- \wedge e^0$ with relations. Moreover, $\Omega^k(A) = 0$ for $k > 3$. In \cite[Theorem 6.3]{AntEmaTho} we show that the following statement holds.
\begin{theorem}
	The coaction $\Delta_{A}:A \to A\otimes H$ extends to $\Omega^{\bullet}(A)$ as a morphism of DGAs. Therefore, $\Omega^{\bullet}(A)$ is a complete differential calculus.
	Moreover, $\Omega^{\bullet}(B)$ is the usual pullback differential calculus on the Podleś sphere.
\end{theorem}
On the Podleś sphere, $\Omega^\bullet(B)$ is induced by restricting to coinvariant forms in $\Omega^1(A)$, namely elements $a e^+ + b e^-$ with $|a| = -2$ and $|b| = 2$. Moreover, in \cite{BeggsMajid} (Proposition 2.35, page 113) it is stated that the volume form can be expressed in terms of elements in $B$ and $\Omega^{1}(B)$. We further elaborate on the corresponding \DJ ur\dj evi\'c braiding in \cite[Section 6.3]{AntEmaTho}.

\paragraph{Smash product algebras:} 
Recall from Example \ref{ex:HopfGalois} $iii.)$ that cleft extensions are particular cases of QPBs, where the cleaving map $j\colon H\to A$, a convolution invertible right $H$-colinear map, guarantees the existence of the inverse of the Galois map. There is a well-known correspondence \cite{DoiTakeuchi} between cleft extensions and so-called crossed product algebras. The latter are generalizations of the tensor product algebra $B\otimes H$, where the multiplication is twisted by a $2$-cocycle with values in $B$ and a "weak" action. Covariant calculi on crossed product algebras are constructed in \cite{AndreaThomas} and it was shown in \cite[Section 6.4]{AntEmaTho} that these crossed product calculi are complete. Here, we discuss complete calculi on the subclass of crossed product algebras given by smash product algebras. The latter correspond to cleft extensions, where the cleaving map is an algebra morphism, in addition.

Let $H$ be a Hopf algebra and $B$ a left $H$-module algebra. The latter means that $B$ is an (associative unital) algebra, together with a left $H$-action $\cdot\colon H\otimes B\to B$, such that
$$
h\cdot(bb')=(h_1\cdot b)(h_2\cdot b'),\qquad\qquad
h\cdot 1_B=\varepsilon(h)1_B
$$  
for all $h\in H$ and $b,b'\in B$. Then, the \textit{smash product algebra} is the vector space $B\otimes H$, endowed with the associative multiplication
\begin{equation}
	(b\otimes h)\cdot_\#(b'\otimes h')
	:=b(h_1\cdot b')\otimes h_2h',
\end{equation}
the so-called \textit{smash product}, and unit $1\otimes 1_B$. We denote $B\#H:=(B\otimes H,\cdot_\#)$ and often write $b\# h$ for elements of $B\#H$.

Following \cite{PflaumSchauenburg} we construct a right $H$-covariant differential calculus on $B\#H$ from a bicovariant FODC $(\Omega^1(H),\mathrm{d}_H)$ on $H$ with maximal prolongation $\Omega^\bullet(H)$ and a left $H$-module differential calculus $\Omega^\bullet(B)$ on $B$. The latter is a DC with left $H$-actions $\cdot\colon H\otimes\Omega^k(B)\to\Omega^k(B)$ for all $k\geq 0$, such that 
$$
h\cdot(\omega\wedge\eta)=(h_1\cdot\omega)\wedge(h_2\cdot\eta),\qquad\qquad\text{for all }h\in H,\quad\omega,\eta\in\Omega^\bullet(B)
$$
and such that $\mathrm{d}\colon\Omega^k(B)\to\Omega^{k+1}(B)$ is left $H$-linear for all $k\geq 0$. Then $\Omega^\bullet(B\#H):=\bigoplus_{n\geq 0}\Omega^n(B\#H)$ with $\Omega^n(B\#H):=\bigoplus_{k=0}^n\Omega^k(B)\otimes\Omega^{n-k}(H)$is a differential calculus on $B\#H$ with wedge product and differential given by
\begin{align}
	(\omega^B\#\omega^H)\wedge_\#(\eta^B\#\eta^H)
	&:=(-1)^{|\omega^H||\eta^B|}\omega^B\wedge(\omega^H_{-1}\cdot\eta^B)\#\omega^H_0\wedge\eta^H,\\
	\mathrm{d}_\#(\omega^B\#\omega^H)
	&:=\mathrm{d}_B\omega^B\#\omega^H
	+(-1)^{|\omega^B|}\omega^B\#\mathrm{d}_H\omega^H
\end{align}
for all $\omega^B,\eta^B\in\Omega^\bullet(B)$ and $\omega^H,\eta^H\in\Omega^\bullet(H)$. Moreover, $\Omega^\bullet(B\#H)$ becomes right $H$-covariant if endowed with the right $H$-coaction $\Delta_{\Omega^\bullet(B\#H)}:=\mathrm{id}\otimes\Delta_{\Omega^\bullet(H)}\colon\Omega^\bullet(B\#H)\to\Omega^\bullet(B\#H)\otimes H$. We call it the \textit{smash product calculus}. From \cite[Section 6.4]{AntEmaTho} we deduce the following result.
\begin{proposition}
Let $B$ be a left $H$-module algebra, $\Omega^\bullet(H)$ the maximal prolongation of a bicovariant FODC on $H$ and $\Omega^\bullet(B)$ a left $H$-module differential calculus, as above. Then, the smash product calculus $\Omega^\bullet(B\#H)$ is complete. Moreover,
\begin{enumerate}
\item[i.)] basic forms coincide with $\Omega^\bullet(B)$, while horizontal and vertical forms are given by $\mathrm{hor}^\bullet=\Omega^\bullet(B)\otimes H$ and $\mathrm{ver}^\bullet=B\otimes\Omega^\bullet(H)$, respectively.

\item[ii.)] the extended coaction $\Delta_{B\#H}^\bullet\colon\Omega^\bullet(B\#H)\to\Omega^\bullet(B\#H)\otimes\Omega^\bullet(H)$ reads $\Delta_{B\#H}^\bullet=\mathrm{id}\otimes\Delta^\bullet$, where $\Delta^\bullet$ is the extended comultiplication (see Lemma \ref{lem:gradedHopf}).

\item[iii.)] the extended Galois map reads
\begin{equation}
	\tau^\bullet(\omega^H)=(1\#S^\bullet(\omega^H_{[1]}))\otimes_{\Omega^\bullet(B)}(1\#\omega^H_{[2]})
\end{equation}
for all $\omega^H\in\Omega^\bullet(H)$, where $S^\bullet$ is the extended antipode (see Lemma \ref{lem:gradedHopf}).

\item[iv.)] the extended \DJ ur\dj evi\'c braiding reads
\begin{equation*}
\begin{split}
	\sigma^\bullet&((\omega^B\#\omega^H)\otimes_{\Omega^\bullet(B)}(\eta^B\#\eta^H))\\
	&=(-1)^{(|\eta^B|+|\eta^H|)(|\omega^H_{[2]}|+|\omega^H_{[3]}|)}(\omega^B\#\omega^H_{[1]})\wedge_\#(\eta^B\#\eta^H\wedge S^\bullet(\omega^H_{[2]}))\otimes_{\Omega^\bullet(B)}(1\#\omega^H_{[3]})\\
	&=(-1)^{|\eta^B||\omega^H|+|\eta^H|(|\omega^H_{[2]}|+|\omega^H_{[3]}|)}\\
	&\qquad\qquad(\omega^B\wedge((\omega^H_{[1]})_{-1}\cdot\eta^B)\#(\omega^H_{[1]})_0\wedge\eta^H\wedge S^\bullet(\omega^H_{[2]}))\otimes_{\Omega^\bullet(B)}(1\#\omega^H_{[3]})
\end{split}
\end{equation*}
for all $\omega^B,\eta^B\in\Omega^\bullet(B)$ and $\omega^H,\eta^H\in\Omega^\bullet(H)$.
\end{enumerate}
\end{proposition}
As mentioned before, there is a generalization of the above proposition to crossed product algebras discussed in \cite[Section 6.4]{AntEmaTho}. Examples of smash product calculi can be found in \cite{PflaumSchauenburg} and example of crossed product calculi in \cite{AndreaThomas}.

\section{Gauge transformations}

In the context of quantum principal bundles there is a notion of quantum gauge transformation proposed in \cite{Brzezinksi}. The idea is to understand gauge transformations of a QPB $B=A^{\mathrm{co}H}\subseteq A$ as unital, convolution invertible, colinear maps $H\to A$. As in the classical case, quantum gauge transformations form a group, which turns out to be isomorphic to the group of vertical automorphisms, i.e. unital, left $B$-linear and right $H$-colinear bijections $A\to A$. Following an idea proposed in \cite{Amilcar}, we extend this picture to differential forms for complete calculi. This then allows to act with quantum gauge transformations on connections and their curvatures. If this extension is compatible with the DGA structure, the transformed curvature corresponds to the curvature of the transformed connection. At the end we discuss quantum gauge transformations and connections on the noncommutative $2$-torus.

\subsection{Gauge transformations on quantum principal bundles}

We introduce the concept of gauge transformation, following \cite[Section 5]{Brzezinksi}. Let $B:=A^{\mathrm{co}H}\subseteq A$ be a quantum principal bundle.
\begin{definition}\label{def:gau}
A \textbf{gauge transformation} of the quantum principal bundle $B\subseteq A$ is a $\Bbbk$-linear map $f\colon H\to A$, such that
\begin{enumerate}
\item[i.)] $f$ is convolution invertible, i.e. there exists a $\Bbbk$-linear map $f^{-1}\colon H\to A$, such that 
$$
f(h_1)f^{-1}(h_2)=\varepsilon(h)1_A=f^{-1}(h_1)f(h_2)
$$ 
for all $h\in H$.

\item[ii.)] $f$ is unital, i.e. $f(1)=1_A$.

\item[iii.)] $f$ is right $H$-colinear, where we endow $H$ with the adjoint right coaction $\mathrm{Ad}\colon H\to H\otimes H$, $\mathrm{Ad}(h)=h_2\otimes S(h_1)h_3$, i.e.
\begin{equation}
	\Delta_A\circ f=(f\otimes\mathrm{id})\circ\mathrm{Ad}
\end{equation}
holds.
\end{enumerate}
We denote the set of gauge transformations of $B\subseteq A$ by $\mathfrak{Gau}(B,A)$.
\end{definition}
Note that the convolution product $*$ gives $\mathfrak{Gau}(B,A)$ a group structure. In fact, given two gauge transformations $f,g\in\mathfrak{Gau}(B,A)$, their convolution product $f*g\colon H\to A$, $(f*g)(h):=f(h_1)g(h_2)$ is a gauge transformation. Its convolution inverse is given by $g^{-1}*f^{-1}$, it is clearly unital, and
\begin{align*}
	(f*g\otimes\mathrm{id})(\mathrm{Ad}(h))
	&=(f*g)(h_2)\otimes S(h_1)h_3\\
	&=f(h_2)g(h_3)\otimes S(h_1)h_4\\
	&=f(h_2)g(h_5)\otimes S(h_1)h_3S(h_4)h_6\\
	&=(f(h_2)\otimes S(h_1)h_3)(g(h_5)\otimes S(h_4)h_6)\\
	&=\Delta_A(f(h_1))\Delta_A(g(h_2))\\
	&=\Delta_A(f(h_1)g(h_2))\\
	&=\Delta_A(f*g(h))
\end{align*}
shows that $f*g$ is also right $H$-colinear. The inverse of $f\in\mathfrak{Gau}(B,A)$ is its convolution inverse $f^{-1}$ and the unit of $\mathfrak{Gau}(B,A)$ is the convolution unit $\eta_A\circ\varepsilon_H\colon H\to A$, $h\mapsto\varepsilon(h)1_A$.

It turns out that gauge transformations correspond to vertical automorphisms of the quantum principal bundle.
\begin{definition}
A \textbf{vertical automorphism} of $B\subseteq A$ is a $\Bbbk$-linear map $F\colon A\to A$, such that
\begin{enumerate}
\item[i.)] $F$ is bijective and left $B$-linear.

\item[ii.)] $F$ is unital, i.e. $F(1)=1$.

\item[iii.)] $F$ is right $H$-colinear, i.e.
\begin{equation}
	\Delta_A\circ F=(F\otimes\mathrm{id})\circ\Delta_A
\end{equation}
holds.
\end{enumerate}
We denote the set of vertical automorphisms of $B\subseteq A$ by $\mathfrak{aut}_v(B,A)$.
\end{definition}
Note that $\mathfrak{aut}_v(B,A)$ is a group with respect to the opposite composition $F\otimes G\mapsto G\circ F$ and group inverse given by the $\Bbbk$-linear inverse.
\begin{proposition}\label{prop:gauge}
There is a group isomorphism $\mathfrak{Gau}(B,A)\cong\mathfrak{aut}_v(B,A)$ given by
\begin{equation}
\begin{split}
	\theta\colon\mathfrak{Gau}(B,A)\to\mathfrak{aut}_v(B,A),\qquad
	\theta(f):=F_f\colon &A\to A,\\
	&a\mapsto F_f(a):=a_0f(a_1).
\end{split}
\end{equation}
Its inverse is
\begin{equation}
\begin{split}
\theta^{-1}\colon\mathfrak{aut}_v(B,A)\to\mathfrak{Gau}(B,A),\qquad
\theta^{-1}(F):=f_F\colon &H\to A,\\
&h\mapsto h^{\langle 1\rangle}F(h^{\langle 2\rangle}).
\end{split}
\end{equation}
\end{proposition}
\begin{proof}
For $f\in\mathfrak{Gau}(B,A)$ we obtain an element $\theta(f)=F_f\in\mathfrak{aut}_v(B,A)$. In fact, $F_f$ is clearly left $B$-linear and unital. Its $\Bbbk$-linear inverse is given by $(F_f)^{-1}(a):=a_0f^{-1}(a_1)$, where $f^{-1}\colon H\to A$ denotes the convolution inverse of $f$. The latter is the case, since
\begin{align*}
	(F_f)^{-1}(F_f(a))
	&=(F_f)^{-1}(a_0f(a_1))
	=a_0f(a_2)_0f^{-1}(a_1f(a_2)_1)
	=a_0f(a_3)f^{-1}(a_1S(a_2)a_4)\\
	&=a_0f(a_1)f^{-1}(a_2)
	=a
\end{align*} 
and
\begin{align*}
	F_f((F_f)^{-1}(a))
	&=F_f(a_0f^{-1}(a_1))
	=a_0f^{-1}(a_2)_0f(a_1f^{-1}(a_2)_1)\\
	&=a_0f^{-1}(a_3)f(a_1S(a_2)a_4)
	=a.
\end{align*}
Lastly, $F_f$ is right $H$-colinear, which follows from
\begin{align*}
	\Delta_A(F_f(a))
	&=\Delta_A(a_0f(a_1))
	=a_0f(a_2)_0\otimes a_1f(a_2)_1
	=a_0f(a_3)\otimes a_1S(a_2)a_4
	=af(a_1)\otimes a_2\\
	&=(F_f\otimes\mathrm{id})\Delta_A(a),
\end{align*}
and thus $F_f$ is a vertical automorphism.

Conversely, given a vertical automorphism $F\in\mathfrak{aut}_v(B,A)$ we obtain a gauge transformation $f_F\in\mathfrak{Gau}(B,A)$. First of all, note that $f_F(h)=h^{\langle 1\rangle}F(h^{\langle 2\rangle})$ is well-defined as the composition $f_F=m\circ(\mathrm{id}_A\otimes_BF)\circ\tau$ since $F$ is left $B$-linear. The convolution inverse of $f_F$ is given by
\begin{equation}
	f_F^{-1}\colon H\to A,\qquad\qquad
	f_F^{-1}(h):=h^{\langle 1\rangle}F^{-1}(h^{\langle 2\rangle}),
\end{equation}
where $F^{-1}\colon A\to A$ is the $\Bbbk$-linear inverse of $F$. In fact,
\begin{align*}
	f_F(h_1)f_F^{-1}(h_2)
	&=h_1{}^{\langle 1\rangle}\overbrace{F(h_1{}^{\langle 2\rangle})h_2{}^{\langle 1\rangle}}^{\in B}F^{-1}(h_2{}^{\langle 2\rangle})\\
	&=h_1{}^{\langle 1\rangle}F^{-1}(F(h_1{}^{\langle 2\rangle})h_2{}^{\langle 1\rangle}h_2{}^{\langle 2\rangle})\\
	&=h^{\langle 1\rangle}F^{-1}(F(h^{\langle 2\rangle}))\\
	&=h^{\langle 1\rangle}h^{\langle 2\rangle}\\
	&=\varepsilon(h)1,
\end{align*}
using the properties of the translation map \eqref{translation} and that $F^{-1}$ is left $B$-linear. Similarly $f_F^{-1}(h_1)f_F(h_2)=\varepsilon(h)1$ follows.

Moreover,
$$
\theta(\theta^{-1}(F))(a)
=a_0f_F(a_1)
=a_0(a_1)^{\langle 1\rangle}F((a_1)^{\langle 2\rangle})
=F(a)
$$
and
$$
\theta^{-1}(\theta(f))(h)
=h^{\langle 1\rangle} F_f(h^{\langle 2\rangle})
=h^{\langle 1\rangle}(h^{\langle 2\rangle})_0f((h^{\langle 2\rangle})_1)
=f(h)
$$
show that $\theta$ is a bijection and
$$
\theta(f*g)(a)
=a_0(f*g)(a_1)
=a_0f(a_1)g(a_2)
=\theta(g)(\theta(f)(a))
=(\theta(f)\theta(g))(a)
$$
for all $f,g\in\mathfrak{Gau}(B,A)$, $F\in\mathfrak{aut}_v(B,A)$, $a\in A$ and $h\in H$, proves that $\theta$ is a group morphism.
This concludes the proof.
\end{proof}

\subsection{Extension of gauge transformations to forms}

In this section we follow the idea of \cite[Section 4]{Amilcar} to extend the concepts of gauge transformations and vertical automorphisms introduced in \cite{Brzezinksi} to differential calculi, using the extended coaction approach \cite{DurII} of \DJ ur\dj evi\'c. A correspondence of gauge transformations and vertical automorphisms in the spirit of the previous section will be given.

Recall that a $\Bbbk$-linear map $f^\bullet\colon V^\bullet\to W^\bullet$ between $\mathbb{N}$-graded $\Bbbk$-vector space is of degree zero, if $f^\bullet(V^n)\subseteq W^n$ for all $n\geq 0$. In this case we identify $f^\bullet$ with the collection $\{f^n\}_{n\in\mathbb{N}_0}$ of maps $f^n\colon V^n\to W^n$.

Fix a complete calculus $\Omega^\bullet(A)$ on a QPB $B:=A^{\mathrm{co}H}\subseteq A$.
\begin{definition}\label{def:gradedGauge}
We define $\mathfrak{Gau}^\bullet(B,A)$ as the $\Bbbk$-linear maps $f^\bullet\colon\Omega^\bullet(H)\to\Omega^\bullet(A)$ of degree zero, such that
\begin{enumerate}
\item[i.)] $f^\bullet$ is convolution invertible, i.e. there exists a $\Bbbk$-linear map $\overline{f}^\bullet\colon\Omega^\bullet(H)\to\Omega^\bullet(A)$ of degree zero, such that 
$$
f^\bullet(\omega_{[1]})\overline{f}^\bullet(\omega_{[2]})=\varepsilon^\bullet(\omega)1_A=\overline{f}^\bullet(\omega_{[1]})f^\bullet(\omega_{[2]})
$$ 
for all $\omega\in\Omega^\bullet(H)$.

\item[ii.)] $f^0(1)=1_A$.

\item[iii.)] $\Delta_A^\bullet\circ f^\bullet
=(f^\bullet\otimes\mathrm{id})\circ\mathrm{Ad}^\bullet$, 
where $\mathrm{Ad}^\bullet\colon\Omega^\bullet(H)\to\Omega^\bullet(H)\otimes\Omega^\bullet(H)$, denotes the adjoint coaction $\mathrm{Ad}^\bullet(\omega):=(-1)^{|\omega_{[1]}||\omega_{[2]}|}\omega_{[2]}\otimes S^\bullet(\omega_{[1]})\wedge\omega_{[3]}$.
\end{enumerate}
We call $\mathfrak{Gau}^\bullet(B,A)$ the \textbf{graded gauge transformations}.
\end{definition}
$\mathfrak{Gau}^\bullet(B,A)$ forms a ($\mathbb{N}$-graded) group with respect to the convolution product and convolution unit. It turns out that condition $i.)$ of Definition \ref{def:gradedGauge} automatically extends to higher degrees.
\begin{proposition}\label{prop:degree0}
	Let $f^\bullet\colon\Omega^\bullet(H)\to\Omega^\bullet(A)$ be a $\Bbbk$-linear map of degree zero. Then $f^\bullet=\bigoplus\nolimits_{n\geq 0}f^n$ is convolution invertible if and only if $f^0\colon H\to A$ is convolution invertible. In particular, a degree zero map $f^\bullet\colon\Omega^\bullet(H)\to\Omega^\bullet(A)$ is a graded gauge transformation if and only if $f^0\colon H\to A$ is a gauge transformation and Definition \ref{def:gradedGauge} $iii.)$ is satisfied. 
\end{proposition}
\begin{proof}
Let $f^\bullet=\bigoplus\nolimits_{n\geq 0}f^n\colon\Omega^\bullet(H)\to\Omega^\bullet(A)$ be a degree $0$ map, i.e. $f^n\colon\Omega^n(H)\to\Omega^n(A)$, and assume that $f^0\colon H\to A$ is convolution invertible with convolution inverse $f^{-1}\colon H\to A$. Assume that $g^\bullet=\bigoplus\nolimits_{n\geq 0}g^n\colon\Omega^\bullet(H)\to\Omega^\bullet(A)$ is the convolution inverse of $f^\bullet$. In particular $g^0=f^{-1}$. Furthermore,
$$
0
=\varepsilon(\omega)1_A
=f^\bullet(\omega_{[1]})g^\bullet(\omega_{[2]})
=f^0(\omega_{-1})g^1(\omega_0)+f^1(\omega_0)g^0(\omega_1)
$$
for $\omega\in\Omega^1(H)$ implies that
$$
g^1(\omega)=-f^{-1}(\omega_{-1})f^1(\omega_0)f^{-1}(\omega_1).
$$
For $\omega=h\mathrm{d}h'\wedge\mathrm{d}h''\in\Omega^2(H)$ we obtain
\begin{align*}
	0
	&=\varepsilon(\omega)1_A
	=f^\bullet(\omega_{[1]})g^\bullet(\omega_{[2]})\\
	&=f^0(h_1h'_1h''_1)g^2(h_2\mathrm{d}h'_2\wedge\mathrm{d}h''_2)
	+f^2(h_1\mathrm{d}h'_1\wedge\mathrm{d}h''_1)f^{-1}(h_2h'_2h''_2)\\
	&\quad+f^1(h_1\mathrm{d}(h'_1)h''_1)g^1(h_2h'_2\mathrm{d}h''_2)
	-f^1(h_1h'_1\mathrm{d}h''_1)g^1(h_2\mathrm{d}(h'_2)h''_2)
\end{align*}
and thus we can express the second order $g^2$ of $g^\bullet$ as
\begin{align*}
	g^2(h\mathrm{d}h'\wedge\mathrm{d}h'')
	&=-f^{-1}(h_1h'_1h''_1)f^2(h_2\mathrm{d}h'_2\wedge\mathrm{d}h''_2)f^{-1}(h_3h'_3h''_3)\\
	&\quad-f^{-1}(h_1h'_1h''_1)f^1(h_2\mathrm{d}(h'_2)h''_2)g^1(h_3h'_3\mathrm{d}h''_3)\\
	&\quad+f^{-1}(h_1h'_1h''_1)f^1(h_2h'_2\mathrm{d}h''_2)g^1(h_3\mathrm{d}(h'_3)h''_3).
\end{align*}
Similarly one can express the higher orders of $g^\bullet$ inductively in terms of $f^\bullet$ and $f^{-1}$. Thus, we can construct the convolution inverse of $f^\bullet$ inductively from $f^\bullet$ and $f^{-1}$. This concludes the proof.
\end{proof}
We continue by defining the graded analogue of vertical automorphisms.
\begin{definition}
We further define $\mathfrak{aut}_v^\bullet(B,A)$ as the $\Bbbk$-linear maps $F^\bullet\colon\Omega^\bullet(A)\to\Omega^\bullet(A)$ of degree zero, such that
\begin{enumerate}
\item[i.)] $F^\bullet$ is left $\Omega^\bullet(B)$-linear and bijective.

\item[ii.)] $F^\bullet$ is unital, i.e. $F^0(1)=1$.

\item[iii.)] $F^\bullet$ is right $\Omega^\bullet(H)$-colinear, i.e. $\Delta_A^\bullet\circ F^\bullet=(F^\bullet\otimes\mathrm{id})\circ\Delta_A^\bullet$.
\end{enumerate}
We call $\mathfrak{aut}_v^\bullet(B,A)$ the \textbf{graded vertical automorphisms}.
\end{definition}
There is a ($\mathbb{N}$-graded) group structure on $\mathfrak{aut}_v^\bullet(B,A)$ given by the opposite composition $F^\bullet\otimes G^\bullet\mapsto G^\bullet\circ F^\bullet$ and the $\Bbbk$-linear inverse.

As in the previous section, graded gauge transformations and graded vertical automorphisms are in bijective correspondence.
\begin{proposition}\label{prop:gradedgauge}
There is an isomorphism of ($\mathbb{N}$-graded) groups $\mathfrak{Gau}^\bullet(B,A)\cong\mathfrak{aut}_v^\bullet(B,A)$. Explicitly,
\begin{equation}
	\begin{split}
		\theta^\bullet\colon\mathfrak{Gau}^\bullet(B,A)\to\mathfrak{aut}_v^\bullet(B,A),\qquad
		\theta^\bullet(f^\bullet):=F^\bullet_{f^\bullet}\colon &\Omega^\bullet\to\Omega^\bullet(A),\\
		&a\mapsto F^\bullet_{f^\bullet}(\omega):=\omega_{[0]}f^\bullet(\omega_{[1]})
	\end{split}
\end{equation}
with inverse
\begin{equation}
	\begin{split}
		\theta^{\bullet-1}\colon\mathfrak{aut}_v^\bullet(B,A)\to\mathfrak{Gau}^\bullet(B,A),\qquad
		\theta^{\bullet-1}(F^\bullet):=f^\bullet_{F^\bullet}\colon &\Omega^\bullet(H)\to\Omega^\bullet(A),\\
		&\omega\mapsto\omega^{\langle 1\rangle}F^\bullet(\omega^{\langle 2\rangle}).
	\end{split}
\end{equation}
\end{proposition}
\begin{proof}
This follows in complete analogy to the proof of Proposition \ref{prop:gauge}.
\end{proof}

In the following we will be mainly interested in graded gauge transformations and graded vertical automorphisms up to degree $1$. According to Proposition \ref{prop:degree0} a map of degree zero $f^\bullet=f^0\oplus f^1\colon\Omega^{\leq 1}(H)\to\Omega^{\leq 1}(A)$ is a graded gauge transformation if and only if $f^0\colon H\to A$ is a gauge transformation and
\begin{equation*}
\begin{split}
	\Delta_{\Omega^1(A)}(f^1(\omega))
	&=f^1(\omega_0)\otimes S(\omega_{-1})\omega_1\\
	\mathrm{ver}^{0,1}(f^1(\omega))
	&=f^0(\omega_{-1})\otimes S(\omega_{-2})\omega_0
	+f^0(\omega_1)\otimes S^\bullet(\omega_0)\omega_2\\
	&=f^0(\omega_{-1})\otimes S(\omega_{-2})\omega_0
	-f^0(\omega_2)\otimes S(\omega_{-1})\omega_0S(\omega_1)\omega_3
\end{split}\qquad
\begin{split}
	&\in\Omega^1(A)\otimes H,\\
	&~\\
	&\in A\otimes\Omega^1(H)
\end{split}
\end{equation*}
hold for all $\omega\in\Omega^1(H)$.

\subsection{Connections and their gauge transformations}

In this section we discuss connections and the action of gauge transformations on them. Recall from Proposition \ref{prop:Atiyah} that, given $B:=A^{\mathrm{co}H}\subseteq A$ a QPB and $\Omega^\bullet(A)$ a complete calculus, then the Atiyah sequence
$$
0\to\mathrm{hor}^1\hookrightarrow\Omega^1(A)\overset{\pi_v}{\twoheadrightarrow}\mathrm{ver}^1\to 0
$$
is exact. A \textit{connection $1$-form} is a right $H$-colinear map $s\colon\Lambda^1\to\Omega^1(A)$, which corresponds to a splitting of the Atiyah sequence. Namely, we demand
\begin{enumerate}
\item[$\bullet$] $\Delta_{\Omega^1(A)}(s(\varpi(h)))=s(\varpi(\pi_\varepsilon(h_2)))\otimes S(h_1)h_3$

\item[$\bullet$] $\pi_v(s(\vartheta))=1_A\otimes\vartheta$
\end{enumerate}
for all $h\in H^+$ and $\vartheta\in\Lambda^1$. We denote the convex set of connection $1$-forms on the QPB $B\subseteq A$ by $\mathrm{con}(B,A)$. 
\begin{lemma}
	Let $B\subseteq A$ be a quantum principal bundle equipped with a first order complete differential calculus. The space of all connections $1$-forms on $B\subseteq A$ is a convex set.  
	\proof 
	Let $s,s':\Lambda^{1}\to \Omega^{1}(A)$ be any two connection $1$-forms on $B\subseteq A$, and let $t\in\Bbbk$. We show that the convex combination $\tilde{s}:=ts+(1-t)s'\colon\Lambda^{1}\to \Omega^{1}(A)$ is a connection $1$-form. Right $H-$colinearity follows immediately since all maps are linear, therefore we are left to check that $\pi_{v}\circ \tilde{s} = 1_{A} \otimes \mathrm{id}$. Let $h\in H^{+}:=\ker\varepsilon$. We find 
	\begin{equation}\nonumber 
		\begin{aligned}
			\pi_{v}\circ \tilde{s}\, (\varpi(h)) & = t \, \pi_{v} \circ s\, (\varpi(h))+ (1-t)\, \pi_{v}\circ s'(\varpi(h)) \\ 
			& = t\otimes \varpi(h) + (1-t)\otimes \varpi(h) \\ 
			& = 1\otimes \varpi(h)\, , 
		\end{aligned}
	\end{equation} 
	that is $\pi_{v}\circ \tilde{s} = 1_{A}\otimes \mathrm{id}$. \qed 
\end{lemma}
Given a connection $1$-form $s\colon\Lambda^1\to\Omega^1(A)$ we obtain a \textit{connection} $\Pi\colon\Omega^1(A)\to\Omega^1(A)$ by
\begin{equation}
	\Pi(a\mathrm{d}a'):=aa'_0s(\varphi(\pi_\varepsilon(a'_1)))
\end{equation}
for all $a\mathrm{d}a'\in\Omega^1(A)$. Recall that $\Pi$ is a morphism in ${}_A\mathcal{M}^H$, satisfying
\begin{enumerate}
\item[$\bullet$] $\Pi^2=\Pi$

\item[$\bullet$] $\ker\Pi=\mathrm{hor}^1$.
\end{enumerate}
A connection $\Pi\colon\Omega^1(A)\to\Omega^1(A)$ is called \textit{strong}, if
\begin{equation}
	(\mathrm{id}-\Pi)(\mathrm{d}A)\subseteq\Omega^1(B)A
\end{equation}
holds. The theory of connections and connection $1$-forms is well-known, see e.g. \cite[Chapter 5]{BeggsMajid}, \cite{DaGrHa} and references therein. The bijective correspondence of connection $1$-forms, connections and splittings of \eqref{Atiyah} is for example discussed in \cite[Proposition 4.3]{AntEmaTho}.

Given a right $H$-comodule $V$ one constructs the associated bundle as the coinvariant subspace $E:=(A\otimes V)^{\mathrm{co}H}$. Following \cite[Proposition 5.48]{BeggsMajid} we obtain for every strong connection $\Pi\colon\Omega^1(A)\to\Omega^1(A)$ a \textit{covariant derivative}
\begin{equation}
	\nabla\colon E\to\Omega^1(B)\otimes_BE,\qquad\qquad
	\nabla(a\otimes v):=(\mathrm{id}-\Pi)(da)\otimes v
\end{equation}
on $E$. Above, we identified $(\Omega^1(B)A\otimes V)^{\mathrm{co}H}\cong\Omega^1(B)\otimes(A\otimes V)^{\mathrm{co}H}$, with one map being the module action and its inverse being
$$
\theta\otimes v\mapsto\theta_0\cdot(\theta_1)^{(1)}\otimes(\theta_1)^{(2)}\otimes v
$$
for all $\theta\otimes v\in(\Omega^1(B)A\otimes V)^{\mathrm{co}H}$, where $h\mapsto h^{(1)}\otimes h^{(2)}\in A\otimes A$ denotes the strong connection obtained from the faithful flatness assumption (see \cite[Section 4.2]{AntEmaTho} for details).
The \textit{curvature} of $\nabla$ is
\begin{equation}
	\mathrm{R}_\nabla(a\otimes v):=-a_0\mathrm{R}_s(\pi_\varepsilon(a_1))\otimes v
\end{equation}
for all $a\otimes v\in(A\otimes V)^{\mathrm{co}H}$, where
\begin{equation}
	\mathrm{R}_s(h):=\mathrm{d}s(\varpi(h))
	+s(\varpi(\pi_\varepsilon(h_1)))\wedge s(\varpi(\pi_\varepsilon(h_2)))
\end{equation}
for all $h\in H^+$ and $s\colon\Lambda^1\to\Omega^1(A)$ is the connection $1$-form associated with $\Pi$.

Inspired by \cite[Section 4]{Amilcar} we define the following action of graded vertical automorphisms $F^\bullet\colon\Omega^\bullet(A)\to\Omega^\bullet(A)$ on connection $1$-forms $s\colon\Lambda^1\to\Omega^1(A)$ by
\begin{equation}
	F^\bullet\rhd s:=F^\bullet\circ s\colon\Lambda^1\to\Omega^1(A).
\end{equation}
\begin{proposition}\label{prop:gaugeact}
$F^\bullet\rhd s$ is a connection $1$-form with curvature
\begin{equation}
	\mathrm{R}_{F^\bullet\rhd s}=\mathrm{d}F^\bullet(s(\varpi(h)))
	+F^\bullet(s(\varpi(\pi_\varepsilon(h_1))))\wedge F^\bullet(s(\varpi(\pi_\varepsilon(h_2)))).
\end{equation}
In particular,
\begin{equation}
	F^\bullet\circ\mathrm{R}_s
	=\mathrm{R}_{F^\bullet\rhd s}
\end{equation}
if $F^\bullet$ is a morphism of DGAs.

Moreover, $\mathfrak{aut}_v(B,A)\otimes\mathrm{con}(B,A)\to\mathrm{con}(B,A)$, $F^\bullet\otimes s\mapsto F^\bullet\rhd s$ is a group action.
\end{proposition}
\begin{proof}
\begin{enumerate}
\item[i.)] \textbf{$F^\bullet\circ s$ is right $H$-colinear:}
recall that
\begin{equation}\label{eq:Fbull}
	\Delta_A^\bullet\circ F^\bullet=(F^\bullet\otimes\mathrm{id})\circ\Delta_A^\bullet\colon\Omega^\bullet(A)\to\Omega^\bullet(A)\otimes\Omega^\bullet(H).
\end{equation}
Projection of the above equation to the component $\Omega^n(A)\otimes H$ (for all $n\geq 0$) shows that $F^n$ is right $H$-colinear. In particular,
\begin{align*}
	\Delta_A(F^1(s(\varpi(h))))
	&=(F^1\otimes\mathrm{id})(\Delta_A(s(\varpi(h))))\\
	&=(F^1\otimes\mathrm{id})(s(\varpi(\pi_\varepsilon(h_2)))\otimes S(h_1)h_3)\\
	&=F^1(s(\varpi(\pi_\varepsilon(h_2))))\otimes S(h_1)h_3
\end{align*}
for all $h\in H^+$.

\item[ii.)] \textbf{$F^\bullet\circ s$ corresponds to a splitting of the Atiyah sequence:} projecting \eqref{eq:Fbull} to $A\otimes\Omega^1(H)$ gives
\begin{align*}
	\mathrm{ver}^{0,1}\circ F^1=(F^0\otimes\mathrm{id})\circ\mathrm{ver}^{0,1}\colon\Omega^1(A)\to A\otimes\Omega^1(H).
\end{align*}
Using the fact that $\pi_v=\kappa\circ\mathrm{ver}^{0,1}\colon\Omega^1(A)\to A\otimes\Lambda^1$, where $\kappa\colon A\square_H\Omega^1(H)\to A\otimes\Lambda^1$, $a\otimes\omega^H\mapsto a_0\otimes S(a_1)\omega^H$, we obtain
\begin{align*}
	\pi_v\circ F^1
	&=\kappa\circ\mathrm{ver}^{0,1}\circ F^1\\
	&=\kappa\circ(F^0\otimes\mathrm{id})\circ\mathrm{ver}^{0,1}\\
	&=(F^0\otimes\mathrm{id})\circ\kappa\circ\mathrm{ver}^{0,1}\\
	&=(F^0\otimes\mathrm{id})\circ\pi_v,
\end{align*}
where in the second to last equality we used the right $H$-colinearity of $F^0$. Thus,
\begin{align*}
	\pi_v(F^1(s(\vartheta)))
	&=(F^0\otimes\mathrm{id})(\pi_v(s(\vartheta)))\\
	&=(F^0\otimes\mathrm{id})(1\otimes\vartheta)\\
	&=1\otimes\vartheta
\end{align*}
for all $\vartheta\in\Lambda^1$, where we also used that $F^0(1)=1$.
\end{enumerate}
\end{proof}
Using Proposition \ref{prop:gradedgauge} in combination with Proposition \ref{prop:gaugeact} we are able to define a group action
\begin{equation}
\begin{split}
	\mathfrak{gau}^\bullet\otimes\mathrm{con}(B,A)&\to(B,A)\\
	f^\bullet\otimes s&\mapsto F^\bullet_{f^\bullet}\circ s
\end{split}
\end{equation}
of graded gauge transformations on connection $1$-forms. Then, utilizing the correspondence of connection $1$-forms, connections and covariant derivatives on associated bundles, it is straight forward to formulate group actions of graded gauge transformations and graded vertical automorphisms on connections, covariant derivatives and their curvatures.

\subsection{Examples}

\noindent We consider the noncommutative algebraic 2-torus as a paradigmatic example to discuss connections and gauge transformations explicitly. As in Section \ref{sec:examples}, let us fix $A:=\mathcal{O}_{\theta}(\mathbb{T}^{2})=\mathbb{C}[u,v,u^{-1},v^{-1}]\mathbin{/}\langle vu-e^{i\theta}uv\rangle$ and $H:=\mathcal{O}(\mathrm{U}(1))=\mathbb{C}[t,t^{-1}]$. The following proposition extends the preliminary result already discussed in \cite[Theorem 6.1]{AntEmaTho}.
\begin{proposition}\label{prop:contorus}
	All connection $1$-forms $s\colon \Lambda^{1}\rightarrow \Omega^{1}(A)$ on the noncommutative algebraic 2-torus are given by convex combinations of 
	\begin{equation}\label{eq:connections_on_torus}
		\begin{aligned}
			t^{-1}\mathrm{d}t & \mapsto u^{-1}\mathrm{d}u + b\mathrm{d}b',\\ 
			t^{-1}\mathrm{d}t & \mapsto v^{-1}\mathrm{d}v + b\mathrm{d}b',
		\end{aligned}
	\end{equation} for any $b,b'\in B:=A^{coH}$.
	\proof 
	The proof is a straightforward check of the axioms of a connection 1-form $s\colon \Lambda^{1}\to \Omega^{1}(A)$. Without loss of generality, we show that the first assignment in Equation \eqref{eq:connections_on_torus} coherently defines connection 1-forms. We have 
	\begin{equation}\nonumber
		\begin{aligned}
			\Delta_{\Omega^{1}(A)}(u^{-1}\mathrm{d}u + b\mathrm{d}b') & = \Delta_{A}(u^{-1})\Delta_{\Omega^{1}(A)}(\mathrm{d}u) +\Delta_{A}(b)\Delta_{\Omega^{1}(A)}(\mathrm{d}b') \\ 
			& = (u^{-1}\otimes t^{-1}) (\mathrm{d}u\otimes t) + (b\otimes 1)(\mathrm{d}b'\otimes 1) \\ 
			& = u^{-1}\mathrm{d}u\otimes 1 + b\mathrm{d}b'\otimes 1\\ 
			& = (u^{-1}\mathrm{d}u + b\mathrm{d}b')\otimes 1\\ 
			& = (s\otimes \mathrm{id})\circ \mathrm{Ad}(t^{-1}\mathrm{d}t),
		\end{aligned}
	\end{equation}
	and moreover
	\begin{equation}\nonumber 
		\begin{aligned}
			\pi_{v}\circ s(t^{-1}\mathrm{d}t) & = \pi_{v}\circ (u^{-1}\mathrm{d}u + b\mathrm{d}b') \\ 
			& = u^{-1}u\otimes S(t^{-1}t)t^{-1}\mathrm{d}t^{n} \\ 
			& = 1_{A}\otimes t^{-1}\mathrm{d}t. 
		\end{aligned}
	\end{equation} 
	Notice how any assignment of the form $t^{-n}\mathrm{d}t^{n}\mapsto u^{-n}\mathrm{d}u^{n}$ defines the same connection by the Leibniz rule of the differential and $\mathbb{k}-$linearity of the connection $1$-form. A similar calculation proves that the second assignment in Equation \eqref{eq:connections_on_torus} defines connection 1-forms. Finally, any other assignment involving both $u$ and $v$ generators fails to be a connection. Indeed, while any combination of the form $u^{m}\mathrm{d}v^{n}$ or $v^{m}\mathrm{d}u^{n}$ for $n,m\in\mathbb{Z}\mathbin{/}\{0\}$ defines a right $H-$coliner map for $n=m$, the second axiom for a connection 1-form does not hold. 
	\qed 
\end{proposition}
We now discuss gauge transformations for the connection 1-form $s(t^{-1}\mathrm{d}t)=u^{-1}\mathrm{d}u$. We notice that condition $iii)$ of Definition \ref{def:gau} is equivalent to asking $f:H\to A$ to map into the coinvariant elements $B:=A^{\mathrm{co}H}$. Indeed 
\begin{equation}\nonumber 
	\begin{aligned}
		\Delta_{A}\circ f(t) &= (f\otimes \mathrm{id}) \circ \mathrm{Ad}(t) \\ 
		&= (f\otimes \mathrm{id}) \circ (t\otimes S(t)t)\\
		&= f(t)\otimes 1_{H}.
	\end{aligned}
\end{equation}
For simplicity we constrain ourselves with determining gauge transformations that are DGA morphisms. Accordingly, without loss of generality, we study gauge transformations
\begin{equation}\label{eq:gauge_torus}
	\begin{aligned}
		f:H&\to A, \\
		t&\mapsto u^{-n}v^{-n},\\ 
		t^{-1}&\mapsto u^{n}v^{n},
	\end{aligned}
\end{equation}
for some integer $n$. Condition $i)$ of Definition \ref{def:gau} is verified with $f^{-1}:H\to A$ sending $t\mapsto v^{n}u^{n}$, as can be easily deduced via the direct calculation 
\begin{equation}\nonumber 
	\begin{aligned}
		f(t_{1})f^{-1}(t_{2}) &= f(t)f^{-1}(t) \\
		&=\varepsilon(t)1_{A}\\
		&=f^{-1}(t_{1})f(t_{2}).
	\end{aligned}
\end{equation} 

Consider now an element $u^k v^\ell$ of $A$ for some $k,\ell \in \mathbb{Z}$. The gauge transformation $f:H\to A$ induces a vertical automorphism $F\colon A\to A$ of the quantum principal bundle $B\subseteq A$ via 
\begin{equation}
	\begin{aligned}\nonumber
		F(u^k v^\ell)&:= \mu\circ (\mathrm{id}\otimes f)\circ\Delta_{A}(u^{k}v^\ell) \\ 
		& = u^k v^\ell f(t^{k-\ell}).
	\end{aligned}
\end{equation}
After this observation we now closely analyse the gauge transformations in Equation \eqref{eq:gauge_torus}. 
Consider the connection $1$-form $s(t^{-1}\mathrm{d}t)=u^{-1}\mathrm{d}u$. By our assumption, $F:A\rightarrow A$ extends as a morphism $F^{\bullet}:\Omega^{\bullet}(A)\to \Omega^{\bullet}(A)$ of DGAs. We then get 
\begin{equation}\nonumber 
	\begin{aligned}
		F^1(u^{-1}du) &= F(u^{-1})\mathrm{d}F(u)\\ 
		&= u^{-1}u^{n}v^{n}\mathrm{d}(uu^{-n}v^{-n})  \\
		&= u^{n-1}v^{n}\mathrm{d}(u^{1-n}v^{-n})\\
		&= u^{n-1}v^{n}((1-n)u^{-n}\mathrm{d}(u)v^{-n} - n \ u^{1-n}v^{-n-1}\mathrm{d}v)\\
		&=(1-n)u^{n-1}v^{n}u^{-n}\mathrm{d}(u)v^{-n}-n\ u^{n-1}v^{n}u^{1-n}v^{-n-1}\mathrm{d}v \\
		& = e^{in(1-n)\theta}\left((1-n)u^{-1}\mathrm{d}u - n\ v^{-1}\mathrm{d}v\right)\\ 
		& = e^{in(1-n)\theta}\left(u^{-1}\mathrm{d}u - n( u^{-1}\mathrm{d}u +  v^{-1}\mathrm{d}v)\right) \\
		& = e^{in(1-n)\theta}\left(u^{-1}\mathrm{d}u + n\ uv\ \mathrm{d}(v^{-1}u^{-1})\right).
	\end{aligned} 
\end{equation}
Notice that $uv\ \mathrm{d}(v^{-1}u^{-1})\in \Omega^1(B)$ is a closed basic form, since, by a simple computation, we obtain $\mathrm{d}(uv\ \mathrm{d}(v^{-1}u^{-1}))=\mathrm{d}(uv)\wedge \mathrm{d}(v^{-1}u^{-1})=0$.

\begin{remark} It is well-known in the literature that, considering a $U(1)$-principal bundle, gauge transformations act on a connection $1$-form by shifting it by a closed base $1$-form. In the example discussed above we recover, up to a phase factor, this classical result.
\end{remark}

We continue by discussing the curvature of the aforementioned connections. As expected, the noncommutative 2-torus is flat.

\begin{proposition}
	All connection 1-forms on the noncommutative 2-torus are flat.
	
	\proof Without loss of generality we may consider the connection $1$-form given by the assignment $t^{-1}\mathrm{d}t \mapsto u^{-1}\mathrm{d}u + b\mathrm{d}b'$. We write the element $t^{-1}\mathrm{d}t$ spanning the coinvariant forms on the structure Hopf algebra via the quantum Cartan--Maurer form as $\varpi(\pi_\varepsilon(t))=S(t_{1})\mathrm{d}(t_{2}) =t^{-1}\mathrm{d}t$.  We find
	\begin{equation}\nonumber
		\begin{aligned}
			R_{s}(h)& =\mathrm{d}\circ s(\varpi(t)) + s(\varpi(t))\wedge s(\varpi(t)) \\ 
			& = \mathrm{d}(s(t^{-1}\mathrm{d}t)) + s(t^{-1}\mathrm{d}t)\wedge s(t^{-1}\mathrm{d}t) \\ 
			& = \mathrm{d}(u^{-1}\mathrm{d}u+b\mathrm{d}b') + (u^{-1}\mathrm{d}u+b\mathrm{d}b')\wedge (u^{-1}\mathrm{d}u+b\mathrm{d}b') \\
			& = \mathrm{d}b\wedge \mathrm{d}b' + u^{-1}\mathrm{d}u \wedge b \mathrm{d}b' + b\mathrm{d}b'\wedge u^{-1}\mathrm{d}u \\
			& = \mathrm{d}(u^{k}v^{k})\wedge \mathrm{d}(u^{\ell}v^{\ell}) + u^{-1}\mathrm{d}u\wedge (u^{k}v^{k})\mathrm{d}(u^{\ell}v^{\ell}) + u^{k}v^{k}\mathrm{d}(u^{\ell}v^{\ell}) \wedge u^{-1}\mathrm{d}u \\ 
			& = k\ell(u^{k-1}\mathrm{d}(u)v^{k}u^{\ell} \wedge v^{\ell-1}\mathrm{d}v + u^{k}v^{k-1}\mathrm{d}v\wedge  u^{\ell-1}\mathrm{d}(u)v^{\ell}) \\ 
			& \quad + \ell (u^{-1}\mathrm{d}u \wedge u^{k}v^{k}u^{\ell}v^{\ell-1}\mathrm{d}v + u^{k}v^{k}u^{\ell}v^{\ell-1}\mathrm{d}v\wedge u^{-1}\mathrm{d}u)\\ 
			& = k\ell e^{-i(k+\ell-1)\theta} u^{k-1}v^{k}u^{\ell}v^{\ell-1} (\mathrm{d}u\wedge\mathrm{d}v + e^{-i\theta}\mathrm{d}v\wedge \mathrm{d}u)\\
			& \quad + \ell e^{-i(k+\ell -1)\theta } u^{k-1}v^{k}u^{\ell}v^{\ell-1}(e^{-i\theta}\mathrm{d}v\wedge\mathrm{d}u + \mathrm{d}u\wedge\mathrm{d}v)\\ 
			& = 0.
		\end{aligned}
	\end{equation}
	By a similar argument we deduce that the connection $1$-form  $t^{-1}\mathrm{d}t \mapsto v^{-1}\mathrm{d}v + b\mathrm{d}b'$ also has vanishing curvature. One easily verifies that convex combinations of the previous connection $1$-forms are flat, as well. Thus, the result follows from Proposition \ref{prop:contorus}.\qed 
\end{proposition}

\section*{Acknowledgements}

We would like to thank Paolo Aschieri, Chiara Pagani, Mauro Mantegazza and Réamonn Ó Buachalla for valuable discussions.
TW is supported by the GA\v{C}R PIF 24-11324I. This pubblication is based on work supported by the European Cooperation in Science and Technology Action 21109 CaLISTA, www.cost.eu, HORIZON-MSCA-2022-SE-01-01 CaLIGOLA, MSCA-DN CaLiForNIA - 101119552.


\end{document}